\RequirePackage{fix-cm}
\documentclass[smallextended,numbook,runningheads]{svjour3}       % onecolumn (second format)
\smartqed  % flush right qed marks, e.g. at end of proof
\usepackage{graphicx}
\usepackage{amssymb}
\usepackage{amsmath}
\usepackage{amsfonts,amssymb,amsbsy,amsmath}
\usepackage{eqnarray}
\usepackage{diagbox}
\usepackage{mathptmx}
\usepackage{float} 
\usepackage{subfigure}
\usepackage[ruled]{algorithm2e}
\usepackage{booktabs}%table

\begin{document}
\title{ Randomized block Krylov space methods for trace and log-determinant estimators\thanks{The work is supported by the National Natural Science Foundation of China (No. 11671060) and the Natural Science Foundation Project of CQ CSTC (No. cstc2019jcyj-msxmX0267)}}



\author{  Hanyu Li  \and
         Yuanyang Zhu 
}

\institute{Hanyu Li \at
              College of Mathematics and Statistics, Chongqing University,
              Chongqing, 401331, P. R. China \\
              \email{lihy.hy@gmail.com}\\
              \email{hyli@cqu.edu.cn}           %  \\
           \and
          Yuanyang Zhu \at
              College of Mathematics and Statistics, Chongqing University,
              Chongqing, 401331, P. R. China\\
              \email{yuanyangzhu926@gmail.com}\\
              \email{yyzhu@cqu.edu.cn}   
}

\date{Received: date / Accepted: date}

\maketitle

\begin{abstract}
%Insert your abstract here. Include keywords, PACS and mathematical
%subject classification numbers as needed.
We present randomized algorithms based on block Krylov space method for estimating the trace and log-determinant of Hermitian positive semi-definite matrices. Using the properties of Chebyshev polynomial and Gaussian random matrix, we provide the %probabilistic 
error analysis of the proposed estimators and obtain the expectation and concentration error bounds. These bounds improve the corresponding ones given in the literature. Numerical experiments are presented to illustrate the performance of the algorithms and to test the error bounds.
\keywords{Randomized algorithm \and  Krylov space method \and Trace estimator \and  Log-determinant estimator \and Chebyshev polynomial}
\subclass{68W20 \and 68W25 \and 15A18 \and 15A15}%MSC 

\end{abstract}

\section{Introduction}
\label{intro}
Computing the trace and the log-determinant of  Hermitian positive semi-definite matrices finds many applications in various problems such as  inverse problem \cite{haber}, generalized cross validation \cite{golub}, %the machine learning \cite{zhang2008}, 
spatial statistics \cite{zhang2008}, and so on. Naturally, it is a straightforward problem if the matrices are explicitly defined and we can access the individual entries. For example, a standard approach for computing the determinant of Hermitian positive definite matrices is to leverage its LU decomposition or Cholesky decomposition \cite[Section 14.6]{higham}. However, in big data age,  it is difficult or expensive to explicitly access the individual entries or we can only access the matrix through matrix vector products. We will focus on the latter case in this paper. For this case, seeking the estimators with high accuracy for trace and log-determinant will be of great interest.
%It is well-known that computing the high dimentional matrices' trace and log-determinant is a frequently encountered problem in various fields, like inverse problem\cite{haber}, parameter estimation, GCV\cite{golub}, machine learning, and spatial statistics\cite{zhang2008}, etc.
%The computations of trace and determinant are easy when the target matrix is explicit, but when the matrix is implicit, we need to find indirect method to estimate trace and determinant. 

For the trace of Hermitian positive semi-definite matrix $A$, the popular and simple estimator is the Monte Carlo estimator proposed by Hutchinson \cite{hutchinson}:
%The commonly used way to estimate the trace of an implicit matrix is Monte-Carlo method. It gives the estimator like
\begin{equation}\nonumber
{\rm Tr}(A)\approx \frac{1}{N}\sum_{i=1}^{N}c_{i}A^{T}c_{i},
\end{equation}
where $N$ is the sample size, and $c_{i}$ are the %indepently random vectors. The initial method is proposed by Hutchinson\cite{hutchinson}, which uses $N$
independently and identically distributed (i.i.d.) Rademacher random vectors.  Hutchinson \cite{hutchinson} showed that this estimator is unbiased. %smallest variance estimator. 
Later, %on the basis of Hutchinson, scholars propose new estimator 
by replacing the i.i.d. random vectors $c_{i}$ with Gaussian random vectors, random unit vectors, or columns of the identity matrix sampled uniformly, some scholars produced some other unbiased estimators \cite{avron,roosta}, where Avron and Toledo \cite{avron} first considered the number of Monte Carlo samples $N$ with which a relative error $\epsilon$ with probability $1-\delta$ can be achieved and defined an ($\epsilon$, $\delta$) estimator:
\begin{equation}\nonumber
\mathbb{P}\left[\|{\rm Tr}(A)-\frac{1}{N}\sum_{i=1}^{N}c_{i}A^{T}c_{i}\|\leqslant \epsilon{\rm Tr}(A)\right]\geqslant 1-\delta.
\end{equation}
Some lower bounds for $N$ for different choice of the random vectors $c_{i}$ were provided in  \cite{avron}, which were improved by Roosta-Khorasani and Ascher \cite{roosta}. 
%and both articles give the lower bound of the sample size $N$ to ensure a good $\left ( \varepsilon ,\delta  \right )$ error bounds. 
In 2017, Lin \cite{lin} proposed two new trace estimators from the view of the randomized low-rank approximation of the matrix $A$ with order $n$:
\begin{equation}\nonumber
{\rm Tr}(A)\approx {\rm Tr}(A\Omega(\Omega^*A\Omega)^\dagger(A\Omega)^*),
\end{equation}
\begin{equation}\nonumber
\quad{\rm Tr}(A)\approx {\rm Tr}(A\Omega((A\Omega)^*A\Omega)^\dagger((A\Omega)^*A(A\Omega))((A\Omega)^*A\Omega)^\dagger(A\Omega)^*),
\end{equation}
where $\Omega$ is a Gaussian random matrix of size $n\times k$ with $k\ll n$ and, for a matrix $X$, $X^\dagger$ denotes its Moore-Penrose inverse. The author mainly investigated the first estimator and found that the method can be much faster than the Monte Carlo estimator. However, there was no formal error analysis for this estimator provided in \cite{lin}. Later, Saibaba et al. \cite{saibaba} also presented a new trace estimator based on  randomized low-rank approximation and 
% to estimate the trace of smooth functions of Hermitian matrices, but did not give a specific error analysis, in the same year, Saibaba et al.\cite{saibaba} 
provided detailed error analysis to validate the theoretical reliability of the estimator. 

For the log-determinant of Hermitian positive definite matrix $A$, a popular approach is to combine the identity % can take advantage of the equation 
$\log \det A= {\rm Tr}\left ( \log A \right )$ with the Monte Carlo estimators for trace introduced above. With this idea, Barry and Pace \cite{barry} first proposed the Monte Carlo estimator of log-determinant for large sparse matrix.  Later, Zhang and Leithead \cite{zhang} generalized the estimator to general Hermitian positive definite matrix. However, both of the above two papers didn't provide the rigorous error analysis of these estimators. Recently,  Boutsidis et al. \cite{boutsidis} continued the above work and investigated the error analysis in detail based on the results from \cite{avron}. In the above log-determinant estimators %in \cite{boutsidis}, 
\cite{barry,zhang,boutsidis},
the Taylor expansion was used to expand $\log(A)$.  Pace and LeSage \cite{pace} first introduced Chebyshev approximation to approximate $\log(A)$, however, they didn't combine the approximation with Monte Carlo estimators for trace and there was no formal error analysis. %So the method  is only applicable to small or medium matrix. 
These works were done by Han et al. \cite{han2015} and they also generalized the method to trace estimator for matrix function \cite{han2017}. In addition, some scholars also used Cauchy integral formula or spline to expand $\log(A)$ and to devise the estimators for log-determinant \cite{aune,anitescu}. Recently, based on randomized low-rank approximation of matrix, Saibaba et al. \cite{saibaba} proposed a new log-determinant estimator without using Taylor expansion or Chebyshev approximation, and discussed the error analysis for this estimator in detail.

% \cite{barry,zhang,boutsidis}. To reduce computational complexity of $\log A$ for large $A$, $\log A$ can be expanded into  a Chebyshev polynomial\cite{pace,han2015,han2017}, or Cauchy integral formula\cite{aune}. Furthermore, Saibaba\cite{saibaba} came up with a new log-determinant estimator based on randomized low-rank approximation of matrices, $\log \det \left ( I_{n}+A \right )\approx \log \det \left ( I_{l}+Q^{*}AQ \right )$.
For the more accurate trace and log-determinant estimators given in \cite{saibaba}, a main and attractive feature is that they took advantage of randomized subspace iteration algorithm, which has been extensively studied and found many applications \cite{halko,gittens,gu2015}.  In recent years, some scholars found that the randomized block Krylov space methods have more advantage compared randomized subspace iteration algorithm \cite{musco,drineas,yuan}. For example, the former has 
%Randomized block Krylov methods have been widely used in low rank matrix approximation in recent years\cite{musco,drineas,yuan}, and it has 
faster eigenvalue convergence rate when the target matrix has a large eigenvalue gap whose location is known. As a result, the randomized block Krylov space method receives more and more attention from the points of
view of theory and applications \cite{musco,drineas,yuan,feng2018}. Inspired by the advantage of the randomized block Krylov space method and the work of Saibaba et al. \cite{saibaba}, we consider the new estimators for trace and log-determinant of Hermitian positive definite matrices and their error analysis in the present paper. The obtained error bounds for these estimators will be tighter than the corresponding ones given in \cite{saibaba} in most of cases.

%From Saibaba's work\cite{saibaba}, we know randomized subspace iteration can provide more accurate estimators for trace and log-determinant than Monte-Carlo method. In this paper, we will use randomized block Krylov methods to give new estimators for $Tr(A)$ and $\log \det \left ( I_{n}+A \right )$ ,also give better error bounds than the predecessors\cite{saibaba}.

The rest of this paper is organized as follows. Section \ref{sec:2} presents some preliminaries. %, including some prior work we want to compare. 
In Section \ref{sec:3}, we provide the main algorithms and the %deterministic 
error analysis of our estimators. The comparisons between our results with the ones from \cite{saibaba} are also discussed in this section. %the deterministic results with randomized subspace iteration's. 
Section \ref{sec:4} is devoted to %probabilistic analysis of the estimators. In Section \ref{Sec5}, we give some 
numerical experiments to illustrate our new randomized estimators and to test the error bounds. %are convergent faster than the estimators given by randomized subspace iteration. 
Finally,  the concluding remarks of the whole paper are presented. 

\section{Preliminaries}
\label{sec:2}
In this section, we first clarify our assumptions. Then, we review some results on Chebyshev polynomials and the algorithms from \cite{saibaba}. %Finally, two useful lemmas are introduced.
\subsection{Assumptions}
\label{sec:2.1}
Throughout this paper, we assume that $A\in \mathbb{C}^{n \times n}$ is a Hermitian %(or real symmetric) 
positive semi-definite matrix %with $k$ dominant eigenvalues, 
and has %a eigenvalue gap. 
the following eigenvalue decomposition: % of a matrix $A$ is defined as 
\begin{align}\label{201}
A= U\Lambda U^{*},\quad \Lambda={\rm diag}\left ( \lambda _{1}\cdot\cdot\cdot \lambda _{n} \right )\in \mathbb{R}^{n \times n},
\end{align}
where  $U\in \mathbb{C}^{n \times n}$ is unitary, the eigenvalues satisfy $\lambda _{1}\geq \cdot \cdot \cdot \geq \lambda _{n}$, and we assume there is a gap in these eigenvalues: $\lambda_{k}>\lambda_{k+1}$.  As done in \cite{saibaba}, we partition $\Lambda$ and $U$ as follows
%To distinguish the dominant eigenvalues from the sub-dominant ones by partitioning
\begin{align*}%\label{20000}
\Lambda= \begin{pmatrix}
\Lambda _{1} &  \\ 
& \Lambda _{2} 
\end{pmatrix},\quad U= \begin{pmatrix}
U_{1} & U_{2}
\end{pmatrix},
\end{align*}
where $ \Lambda _{1}\in \mathbb{R}^{k\times k}
$, $ \Lambda _{2}\in \mathbb{R}^{(n-k)\times (n-k)}
$, $ U_{1}\in \mathbb{C}^{n\times k}$, and $ U_{2}\in \mathbb{C}^{n\times (n-k)}$. %Here, $p$ is the size of oversampling to balance the reliability and the fasterconvergence of algor
%with $ 0\leq p\leq l-k$.  By $p$, we can create a gap: The largest eigenvalue in $ \Lambda _{3}$ is $\lambda_{l-p+1}$, which is potentially much smaller than $\lambda_k$, the smallest eigenvalue of $ \Lambda _{1}$. When $p=l-k$, the above gap disappears and the case reduces the one considered in \cite{saibaba}. For the latter, we need to 

Given a number $q\geqslant 1$ and a Gaussian random matrix $\Omega\in\mathbb{C}^{n \times l}$ with $k\leqslant l=k+p\ll n$, set 
$$K_q=\left(A\Omega, A^2\Omega,\cdots,A^q\Omega\right).$$
Like \cite{drineas}, we write 
$$\mathcal{K}_q={\rm range}\left(A\Omega, A^2\Omega,\cdots,A^q\Omega\right),$$
call it the block Krylov space in $A$ and $\Omega$, and assume  $\dim\left(\mathcal{K}_{q} \right )=ql$. That is, $K_q$ has full column rank.
%\subsection{Krylov subspace}\label{Sec2.3.2}
%To approximate the dominant subspace of $A\in \mathbb{C}^{n \times n}$, we will generate a random Gaussian matrix
%$\Omega\in\mathbb{C}^{n \times b}$, and construct the Krylov subspace of $A$ on $A\Omega$,
%\begin{align}\label{203}
%\mathcal{K}_{q}\equiv\mathcal{K}_{q}\left(A,A\Omega \right)={\rm range}
%\begin{bmatrix}
%A\Omega  &  A^{2}\Omega & \cdots  & A^{q}\Omega 
%\end{bmatrix}.
%\end{align}
%And assume $\dim\left(\mathcal{K}_{q} \right )=qb$. 
%\subsubsection{Matrix function for the Krylov space}\label{Sec2.3.1}
It is known that the elements of the Krylov subspace $\mathcal{K}_{q}$ can be expressed in terms of matrices $ \phi\left(A\right)\Omega\in \mathbb{C}^{n\times l} $ \cite{li2015,drineas}, where $ \phi\left(\cdot\right) $ is a polynomial of degree $q$. Considering the eigenvalue decomposition of $A$ in \eqref{201}, it is easy to verify that %, but we specially take the $ \phi\left(\cdot\right) $ defined in \eqref{202}. And we can denote elements in $\mathcal{K}_{q}$ by
\begin{align*}
K=\phi\left(A\right)\Omega= U\phi\left(\Lambda \right )U^{*}\Omega ,
\end{align*}
where
\begin{align*}
\phi\left(\Lambda \right)={\rm diag}\begin{pmatrix}
\phi\left(\lambda_{1}\right), &\phi\left(\lambda_{2}\right),  &\cdots,  & \phi\left(\lambda_{n}\right)
\end{pmatrix}=\begin{pmatrix}
\phi(\Lambda _{1}) &  \\ 
& \phi(\Lambda _{2}) 
\end{pmatrix}.
\end{align*}
%Obviously, $ {\rm range}\left(K\right)\subset \mathcal{K}_{q} $.\\
%If $ Q_{q} $ is an orthonormal basis of  $\mathcal{K}_{q}$, $ Q $ is an orthonormal basis  for the column span of $ K $, then $ {\rm range}\left(Q\right)\subseteq  {\rm range}\left(Q_{q}\right) $. 

Like \cite{gu2015,saibaba}, we denote $\widehat{\Omega}=U^{*}\Omega $ and hence
%since $U= \begin{pmatrix}
%U_{1} & U_{2}
%\end{pmatrix}$, we have 
\begin{align*}\widehat{\Omega}= \begin{pmatrix}
U_{1}^{*}\Omega \\ 
U_{2}^{*}\Omega
\end{pmatrix} =\begin{pmatrix}
\widehat{\Omega}_{1}\\ 
\widehat{\Omega}_{2}
\end{pmatrix} ,
\end{align*}
where $ \widehat{\Omega}_{1}=U_{1}^{*}\Omega\in\mathbb{C}^{k\times l}$ and $ \widehat{\Omega}_{2}=U_{2}^{*}\Omega\in\mathbb{C}^{(n-k)\times l}$. %, 0\leq p\leq b-k$, parameter $p$ is to balance the need for oversampling for reliability.\\
We assume that $ {\rm rank}( \widehat{\Omega}_{1})=k $, and hence its Moore-Penrose inverse $\widehat{\Omega}_{1}^{\dagger }$ satisfies $\widehat{\Omega}_{1}\widehat{\Omega}_{1}^{\dagger }=I_{k}$.

\subsection{Chebyshev polynomials}
\label{sec:2.2}
%In our deterministic analysis, we will use some knowledge about Chebyshev polynomials.
The $q$th degree Chebyshev polynomial is recursively defined as follows
\begin{align*}
T_{0}\left (x\right)\equiv 1;\quad
T_{1}\left (x\right)\equiv x;\quad
T_{q}\left (x\right)\equiv 2qT_{q-1}\left(x\right)-T_{q-2}\left(x\right).
\end{align*}
They can also be expressed as
\begin{equation*}
T_{q}\left (x\right)=\begin{cases}
\frac{\left ( x+\sqrt{x^{2}-1} \right )^{q}+\left ( x-\sqrt{x^{2}-1} \right )^{q}}{2},\quad x\geq 1, \\
\cos \left ( q\arccos \left ( x \right ) \right ),\quad -1\leq x\leq 1.
\end{cases}
\end{equation*}

%\subsubsection{Reasonable choice of polynomial}\label{Sec2.2.1}
By Chebyshev polynomials, we 
%In order to facilitate the latter proof, we first 
construct a polynomial $f\left (\cdot \right)$ with degree $q-1$, which will play an important role in our error analysis for algorithms.
\begin{align}\label{202}
f\left(x\right)= \frac{T_{q-1}\left ( \frac{2x-\lambda _{k+1}}{\lambda _{k+1}} \right)}{T_{q-1}\left ( \frac{2\lambda _{k}-\lambda _{k+1}}{\lambda _{k+1}} \right)},
\end{align}
where $\lambda _{k+1}$ and $\lambda _{k}$ are the eigenvalues of $A$. %$ T_{q-1}\left (\cdot \right) $ is a $\left (q-1\right )$th degree Chebyshev polynomial, 
By the properties of Chebyshev polynomials \cite{musco,drineas}, we can check that the polynomial $ f\left(x\right) $ has the following properties:
\begin{enumerate}
	\item  $ f\left ( \lambda _{i} \right )\geq 1 $, when $ 1\leq i\leq k $;
	\item  $|f\left ( \lambda _{i} \right )|\leq T_{q-1}^{-1}\left(\frac{2\lambda _{k}-\lambda _{k+1}}{\lambda _{k+1}}\right)$, when $ k+1\leq i\leq n $.
	%\item  $f\left ( \lambda _{i} \right )> 0$, when $ k+1\leq i\leq l-p $ and $p\neq l-k$.
\end{enumerate}
From the property 1, it follows that  $ f\left(\Lambda _{1}\right) $  is nonsingular, and 
\begin{align}\label{203}
\left \| f^{-1}\left(\Lambda _{1}\right) \right \|_{2}\leq 1.
\end{align}
From the property 2, we have that 
\begin{align}\label{204}
\left \| f\left(\Lambda _{2}\right)\right \|_{2}\leq T_{q-1}^{-1}\left(\frac{2\lambda _{k}-\lambda _{k+1}}{\lambda _{k+1}}\right). 
\end{align}
%The property \textit{3).} ensures the inverse matrix $ f^{-1}\left(\Lambda_{2}\right)$ of  $ f\left(\Lambda _{2}\right) $ exist.\\
%Following we will use $q$th polynomial $ \phi\left(x\right) $ defined as 
%\begin{align}\label{202}
%\phi\left(x\right)= xf\left(x\right) .
%\end{align}

\subsection{Randomized subspace iteration algorithm}
\label{sec:2.3}
The following algorithm was proposed by Saibaba et al. \cite{saibaba}. Based on the algorithm, the authors presented the estimators of trace and log-determinant: ${\rm Tr}(A)\approx {\rm Tr}(T)$ and $\log\det(I+A)\approx \log\det(I+T)$. %Compared with the corresponding results in \cite{hutchinson,avron,roosta,lin,barry,zhang}, these two new estimators 

\begin{algorithm}[H]\label{Al1}
	\caption{Randomized subspace iteration \cite{saibaba} }
	%\LinesNumbered %要求显示行号
	\KwIn{$A\in \mathbb{C}^{n \times n}$: Hermitian positive semi-definite matrix; $k$: target rank; $q$: number of subspace iteration; $\Omega \in\mathbb{C}^{n \times l}$: Gaussian random matrix with $k\leq l=k+p\ll n$.
		\\ 
	}
	\KwOut{$T\in \mathbb{C}^{l\times l}$.}
	\begin{enumerate}
		\item Multiply $Y=A^{q}\Omega $;
		\item Thin QR factorization $Y=QR$;
		\item Compute $T=Q^{*}AQ$.
	\end{enumerate}
\end{algorithm}
From this algorithm, it is easy to find that the information in $A\Omega,A^{2}\Omega,\cdots$,
$A^{q-1}\Omega$ is discarded when computing $Q$. Collecting these information and using them for computing $Q$ is one of the main motivations of this study, and is also an important topic in the field of randomized algorithm \cite{musco,drineas,yuan,feng2018}.
%\subsection{Two useful lemmas}\label{Sec2.5}
%\begin{lemma}\label{lem201}{\rm\cite{gu2015}}
%%	Assume that $X\in\mathbb{C}^{l \times l}$ is nonsingular, and the QR factorizations of $K$ and $KX$ are %and $\Omega $ is full column rank. Let 
%	$K=QR$ and $KX= \widetilde{Q}\widetilde{R} $, respectively. Then
%	\begin{align*}
%	QQ^{*}= \widetilde{Q}\widetilde{Q}^{*}.
%	\end{align*}
%where Q is an orthonormal basis of the column span of K.
%\end{lemma}

%\begin{lemma}\label{lem202}{\rm\cite{}}
%	If $ B\in \mathbb{C}^{m \times n} $ and $ C\in \mathbb{C}^{n \times m} $, then
%	\begin{align*}
%	\det(I_{m}\pm BC)=\det(I_{n}\pm CB).
%	\end{align*}
%\end{lemma}

\section{Algorithms and error analysis}
\label{sec:3}
In this section, we first introduce our algorithms for estimating the trace and log-determinant, and then present the error bounds of these new estimators and their proofs. %, the comparisons
%between our results and the ones from\cite{saibaba} are also discussed. Then, the proofs of these theorems are provided step by step, which mainly comprise deterministic analysis and probabilistic  analysis.

\subsection{Algorithms}\label{Sec3.1}
%\subsubsection{The Algorithms}\label{Sec3.1.1}
%Our randomized block Krylov algorithm can find a matrix $T$ which dimension $ ql\ll n $, so that ${\rm Tr}\left( T\right )$ is an estimator for ${\rm Tr}\left(A\right )$, and $\log \det \left ( I_{ql}+T \right )$ is an estimator for $\log \det \left ( I_{n}+A \right )$. 

\begin{algorithm}[H]\label{Al2}
	\caption{Randomized block Krylov space method}
	%\LinesNumbered %要求显示行号
	\KwIn{$A\in \mathbb{C}^{n \times n}$: Hermitian positive semi-definite matrix;  $k$: target rank; $q$: number of block Krylov space iterations; $\Omega \in\mathbb{C}^{n \times l}$: Gaussian random matrix with 
		%$l$, block size, 
		$k\leq l=k+p\ll n$.\\
		%, $q\geq 1$
	}
	\KwOut{$T\in \mathbb{C}^{ql\times ql}$.}
	\begin{enumerate}
		\item Multiply and collect $K_{q}=\left(
		A\Omega,    A^{2}\Omega,  \cdots, A^{q}\Omega 
		\right)$;
		\item Thin QR factorization %Compute an orthonormal basis $Q_{q}$ for the column span of $K_{q}$, using e.g. 
		$K_q=Q_{q}R_{q}$;
		\item Compute $T=Q_{q}^{*}AQ_{q}$.
	\end{enumerate}
\end{algorithm}
Comparing with Algorithm \ref{Al1}, we can find that the randomized block Krylov space method collects the information discarded in Algorithm \ref{Al1} and hence will be more accurate in theory. Numerical experiments in Section \ref{sec:4} also confirm this result. Moreover, the complexity of Algorithm \ref{Al2} is only a little higher than that for Algorithm \ref{Al1}. This is because the main parts of  complexity of the two algorithms, i.e., the complexity in step 1, are the same. Only in steps 2 and 3, the complexity of our algorithm increases. The factor $l$ in complexity is increased to $ql$. Since $q$ is a small number, the total complexity doesn't change very much.  In addition, similar to the discussions in \cite{saibaba}, Algorithm \ref{Al2} is only an idealized version and the idealized block Krylov space iteration can be numerically unstable. In practice, we can alternate matrix products and QR factorizations to tackle this problem \cite[Algorithm 5.2]{saad}.

%We assume that Krylov subspace matrix $ K_{q} $ in \textbf{Algorithm 2} has full column rank, which means $ {\rm rank}\left ( K_{q}\right )= ql $. {\rm \textbf{Algorithm 2}} is an simple version of randomized block Krylov algorithm, its starting guess is a random Gaussian matrix, as we all know, the actual implementation of \textbf{Algorithm 2} needs to consider the loss of orthogonality of the Lanczos vectors, we can overcome this problem by  a reorganization of the computation to use the three-
% term recurrence and bidiagonalization, and reorthogonalizations of the Lanczos vectors at each step.

%\begin{remark}
%	Considering that {\rm \textbf{Algorithm 1}} does not take advantage “priori information”, it has a certain disregard for prior iterates, i.e.  $ A\Omega,\cdots,A^{q-1}\Omega $. Thus, we proposed  {\rm \textbf{Algorithm 2}}, which can take advantage of all previously computed matrix-vector products.
%\end{remark}

%%对算法的说明，注意事项，算法复杂度没有讨论。
\subsection{Error bounds}\label{Sec3.1.2}

\begin{theorem}\label{lem308}{\rm(Expectation bounds)}
	Let $T=Q_{q}^{*}AQ_{q}$ be computed by Algorithm \ref{Al2} and furthermore, let $p\geq 2 $. Then 
	\begin{align*}
	0\leq \mathbb{E} \left[ {\rm Tr}(A)-{\rm Tr}(T)\right]\leq \left(1+\frac{\lambda_{k+1}}{\lambda_{k}}T_{q-1}^{-2}\left(\frac{2\lambda_{k}-\lambda_{k+1}}{\lambda_{k+1}}\right)C_{ge}\right){\rm Tr}({\Lambda_{2}})
	\end{align*}
	and
	\begin{eqnarray*}
		0&\leq& \mathbb{E} \left[ \log\det(I_{n}+A)-\log\det(I_{ql}+T)\right]\\
		&\leq& \log\det\left(I_{n-k}+\frac{\lambda_{k+1}}{\lambda_{k}}T_{q-1}^{-2}\left(\frac{2\lambda_{k}-\lambda_{k+1}}{\lambda_{k+1}}\right)C_{ge}{\Lambda_{2}}\right)+\log\det(I_{n-k}+{\Lambda_{2}}),
	\end{eqnarray*}
	where $C_{ge}= \frac{p+1}{p-1}(\mu +\sqrt{2})^{2}(\frac{1}{2\pi(p+1)})^{\frac{1}{p+1}}(\frac{e\sqrt{l}}{p+1})^{2}$ with $ \mu = \sqrt{n-k}+\sqrt{l}$. %,\xi=\frac{\lambda_{k+1}}{\lambda_{k}},  \gamma =\frac{2\lambda_{k}-\lambda_{k+1}}{\lambda_{k+1}}\geq \frac{1}{\xi}.$
\end{theorem}
%\emph{Proof} 

\begin{theorem}\label{them309}{\rm (Concentration bounds)}
	Let $T=Q_{q}^{*}AQ_{q}$ be computed by Algorithm \ref{Al2} and furthermore, let $ p \geq 2 $. If $ 0\leq \delta \leq 1 $, then with probability at least $ 1-\delta $ 
	\begin{align*}
	0\leq  {\rm Tr}(A)-{\rm Tr}(T)\leq \left(1+\frac{\lambda_{k+1}}{\lambda_{k}}T_{q-1}^{-2}\left(\frac{2\lambda_{k}-\lambda_{k+1}}{\lambda_{k+1}}\right)C_{g}\right){\rm Tr}({\Lambda_{2}})
	\end{align*}
	and
	\begin{eqnarray*}
		0&\leq&  \log\det(I_{n}+A)-\log\det(I_{ql}+T)\\
		&\leq& \log\det\left(I_{n-k}+\frac{\lambda_{k+1}}{\lambda_{k}}T_{q-1}^{-2}\left(\frac{2\lambda_{k}-\lambda_{k+1}}{\lambda_{k+1}}\right)C_{g}{\Lambda_{2}}\right)+\log\det(I_{n-k}+{\Lambda_{2}}),
	\end{eqnarray*}
	where $ C_{g}=\left(\mu +\sqrt{2\log \frac{2}{\delta }}\right)^{2} (\frac{2}{\delta })^{\frac{2}{p+1}}(\frac{e\sqrt{l}}{p+1})^{2}$ with $\mu = \sqrt{n-k}+\sqrt{l}$. %,\xi=\frac{\lambda_{k+1}}{\lambda_{k}},  \gamma =\frac{2\lambda_{k}-\lambda_{k+1}}{\lambda_{k+1}}\geq \frac{1}{\xi}.$
\end{theorem}
%\emph{Proof}   

%\begin{theorem}\label{them302}
%	Let $T=Q_{q}^{*}AQ_{q}$ be computed by {\rm  \textbf{Algorithm 2}}. Then
%	\begin{align*}
%	0\leq {\rm Tr}(A)-{\rm Tr}(T)\leq \eta_{1} {\rm Tr}\left( \Lambda _{3}\right )+{\rm Tr}\begin{pmatrix}
%	\Lambda _{2} & \\ 
%	& \Lambda _{3}
%	\end{pmatrix}
%	\end{align*}
%	where $\eta_{1}=min\left \{ T_{q-1}^{-1}(\gamma)\left \| \widehat{\Omega}_{2}\widehat{\Omega}_{1}^{\dagger}\right\|_{2},\frac{\lambda _{l-p+1}}{\lambda _{k}}T_{q-1}^{-2}(\gamma)\left \| \widehat{\Omega}_{2}\widehat{\Omega}_{1}^{\dagger}\right\|_{2}^{2} \right\}$, and $ \gamma =\frac{2\lambda_{l-p}-\lambda_{l-p+1}}{\lambda_{l-p+1}}\geq \frac{\lambda_{l-p}}{\lambda_{l-p+1}}. $
%\end{theorem}
%\emph{Proof}   See Section \ref{Sec3.3.1}

%\begin{theorem}\label{them303}
%	Let $T=Q_{q}^{*}AQ_{q}$ be computed by {\rm \textbf{Algorithm 2}}. Then
%	\begin{align*}
%	0\leq \log\det(I_{n}+A)-\log\det(I_{ql}+T)\leq \log\det(I_{n-l+p}+\eta_{2}\Lambda_{3})+\log\det(I_{n-k}+\begin{pmatrix}
%	\Lambda_{2} & \\ 
%	& \Lambda_{3}
%	\end{pmatrix})
%	\end{align*}
%	where $\eta_{2}=\frac{\lambda _{l-p+1}}{\lambda _{k}}T_{q-1}^{-2}(\gamma)\left \| \widehat{\Omega}_{2}\widehat{\Omega}_{1}^{\dagger}\right\|_{2}^{2}$, and $ \gamma =\frac{2\lambda_{l-p}-\lambda_{l-p+1}}{\lambda_{l-p+1}}\geq \frac{\lambda_{l-p}}{\lambda_{l-p+1}}. $
%\end{theorem}
%\emph{Proof}   See Section \ref{Sec3.3.2}

%\subsection{Comparison of error bounds }\label{Sec3.2}
%%还没想好咋写
%Based on the aforementioned partition of $\Lambda$, i.e., $\Lambda={\rm diag}(\Lambda_1, \widehat{\Lambda}_2)$,  
In \cite{saibaba}, Saibaba et al. presented the following error bounds for estimators of trace and log-determinant produced by Algorithm \ref{Al1}.
%To compare the error bounds with \cite{saibaba}, we introduce some main results in it, and for convenience, we define
%\begin{align*}
%\widehat{\Lambda_{2}}=\begin{pmatrix}
%\Lambda_{2} & \\ 
%& \Lambda_{3}
%\end{pmatrix}.
%\end{align*}
%Fristly, we introduce two main results of \cite{saibaba} in Lemma \ref{lem306} \ref{lem307} as follows.
\begin{theorem}\label{lem306}{\rm(Expectation bounds) \cite{saibaba}}
	Let $T=Q^{*}AQ$ be computed by Algorithm \ref{Al1} and furthermore, let $p\geq 2 $. Then 
	\begin{align*}
	0\leq \mathbb{E} \left[ {\rm Tr}(A)-{\rm Tr}(T)\right]\leq \left(1+\left(\frac{\lambda_{k+1}}{\lambda_{k}}\right)^{2q-1}C_{ge}\right){\rm Tr}({\Lambda_{2}})
	\end{align*}
	and
	\begin{eqnarray*}
	0&\leq&\mathbb{E} \left[ \log\det(I_{n}+A)-\log\det(I_{l}+T)\right]\\
&\leq&\log\det\left(I_{n-k}+\left(\frac{\lambda_{k+1}}{\lambda_{k}}\right)^{2q-1}C_{ge}{\Lambda_{2}}\right)+\log\det(I_{n-k}+{\Lambda_{2}}).
	\end{eqnarray*}
	%where $C_{ge}\equiv \frac{p+1}{p-1}(\mu _{1}+\sqrt{2})^{2}(\frac{1}{2\pi(p+1)})^{\frac{1}{p+1}}(\frac{e\sqrt{k+p}}{p+1})^{2}$ with $ \mu_{1} \equiv \sqrt{n-k}+\sqrt{k+p} ,\xi=\frac{\lambda_{k+1}}{\lambda_{k}}.$
\end{theorem}
%\emph{Proof}   Find in \cite{saibaba}(Secction 4.1.1).

\begin{theorem}\label{lem307}{\rm (Concentration bounds) \cite{saibaba}}
	Let $T=Q^{*}AQ$ be computed by Algorithm \ref{Al1} and furthermore, let $ p \geq 2 $. If $ 0\leq \delta \leq 1 $, then with probability at least $ 1-\delta $
	\begin{align*}
	0\leq   {\rm Tr}(A)-{\rm Tr}(T)\leq \left(1+\left(\frac{\lambda_{k+1}}{\lambda_{k}}\right)^{2q-1}C_{g}\right){\rm Tr}({\Lambda_{2}})
	\end{align*}
	and
	\begin{eqnarray*}
	0&\leq &\log\det(I_{n}+A)-\log\det(I_{l}+T)\\
	 &\leq&\log\det\left(I_{n-k}+\left(\frac{\lambda_{k+1}}{\lambda_{k}}\right)^{2q-1}C_{g}{\Lambda_{2}}\right)+\log\det(I_{n-k}+{\Lambda_{2}}).
	\end{eqnarray*}
	%where  %$\gamma =\frac{2\lambda_{k}-\lambda_{l-p+1}}{\lambda_{l-p+1}}$ and 
	%$ C_{g}= (\frac{e\sqrt{k+p}}{p+1})^{2}(\frac{2}{\delta })^{\frac{2}{p+1}}(\mu_{1} +\sqrt{2\log \frac{2}{\delta }})^{2}$ with $\mu_{1} \equiv \sqrt{n-k}+\sqrt{k+p} ,\xi=\frac{\lambda_{k+1}}{\lambda_{k}}. $.
\end{theorem}

%\emph{Proof}   Find in \cite{saibaba}(Theorem 2 ).
Note that 
$$\frac{2\lambda_{k}-\lambda_{k+1}}{\lambda_{k+1}}\geq \frac{\lambda_{k}}{\lambda_{k+1}}> 1,$$
and  %the properties of Chebyshev polynomials that 
$T_{q-1}(x)$ increases faster than $x^{q-1}$ when $q\geqslant 3$ and $x> 1$. Thus,  the term
 $$T_{q-1}^{-2}\left({(2\lambda_{k}-\lambda_{k+1})}/{\lambda_{k+1}}\right)$$ in our bounds is smaller than $\left({\lambda_{k+1}}/{\lambda_{k}}\right)^{2q-2}$ in the bounds in Theorems \ref{lem306} and \ref{lem307} when $q\geqslant 3$. However, it must be pointed out that the order of the matrix $T$ in the bounds in Theorems \ref{lem306} and \ref{lem307} is $l$ not $ql$. If we set the order to be $ql$, i.e., set $p=ql-k$, the terms $C_{ge} $ and  $C_{g}$ in Theorems \ref{lem306} and \ref{lem307} will become small and hence the bounds will be reduced. In this case,  our bounds % in Theorem \ref{lem308} and \ref{them309} 
can't be always tighter than the corresponding ones in Theorems \ref{lem306} and \ref{lem307} when $q\geqslant 3$. However, numerical experiments show that in most of cases, our bounds are tighter. The following are some simple examples. We set $n=3000, k=30, p=10, {\lambda_{k}}/{\lambda_{k+1}}=20$, and $\delta=0.01$, and $q=3,4,5$, respectively. Upon some computations, we have %Naturally, our bounds in Theorems \ref{them301} and \ref{them3002} have more advantages.
Table \ref{tab:1}, where $C_{geq}$ and $C_{gq}$ are derived from $C_{ge}$ and $C_{g}$, respectively, by replacing  $p$ with $ql-k$.

\begin{table}[h]
		\caption{Comparison of terms appearing in bounds for different values of $q$}
	\label{tab:1}       % Give a unique label
	\begin{tabular}{cccc}
		\toprule
		$q$& 3 & 4 & 5\\
	\noalign{\smallskip}\hline\noalign{\smallskip}
		$\frac{\lambda_{k+1}}{\lambda_{k}}T_{q-1}^{-2}\left(\frac{2\lambda_{k}-\lambda_{k+1}}{\lambda_{k+1}}\right)C_{ge}$ & 4.2541e-05 & 6.9946e-09 & 1.1500e-12\\
		\hline
	$\left(\frac{\lambda_{k+1}}{\lambda_{k}}\right)^{2q-1}C_{geq}$& 1.4266e-04 & 2.4408e-07 & 4.7055e-10\\
	\hline
		$\frac{\lambda_{k+1}}{\lambda_{k}}T_{q-1}^{-2}\left(\frac{2\lambda_{k}-\lambda_{k+1}}{\lambda_{k+1}}\right)C_{g}$& 1.4210e-04 & 2.3363e-08 & 3.8414e-12\\
		\hline
		$	\left(\frac{\lambda_{k+1}}{\lambda_{k}}\right)^{2q-1}C_{gq}$& 1.7747e-04 & 2.8924e-07 & 5.4284e-10\\
		\bottomrule
	\end{tabular}
\end{table}

\subsection{Proofs of Theorems \ref{lem308} and \ref{them309}}
\label{sec:3.3}

We only prove the structural (deterministic) error bounds for trace and log-determinant estimators, i.e., we consider $\Omega$ to be any matrix satisfying assumptions given in Section \ref{sec:2.1}. %, deterministic matrix. 
The final results, i.e., the expectation and concentration error bounds in Theorems \ref{lem308} and \ref{them309}, can be derived immediately as done in \cite[Sections 4.1.1 and 4.1.2]{saibaba} by combining the structural error bounds in Theorems \ref{them310} and \ref{them311} given below and \cite[Lemmas 4 and 5]{saibaba}. We first do some preparation for these proofs. A useful lemma is listed as follows.
\begin{lemma}\label{lem201}{\rm\cite{gu2015}}
	Assume that $X\in\mathbb{C}^{l \times l}$ is nonsingular, and the thin QR factorizations of $K$ and $KX$ are %and $\Omega $ is full column rank. Let 
	$K=QR$ and $KX= \widetilde{Q}\widetilde{R} $, respectively. Then
	\begin{align*}
	QQ^{*}= \widetilde{Q}\widetilde{Q}^{*}.
	\end{align*}
	%where Q is an orthonormal basis of the column span of K.
\end{lemma}

This simple result plays an important role in deriving error bounds because we can choose a special $X$ to achieve the useful information of range of $K$. This technique was proposed by Gu \cite{gu2015}. Now we introduce how to find the special $X$. %we will find a specially constructed, full rank matrix $X$, which satisfy Lemma \ref{lem201}.\\
%Denote \begin{align*} \widehat{\Omega}=U^{*}\Omega ,\end{align*}
%since $U= \begin{pmatrix}
%U_{1} & U_{2}
%\end{pmatrix}$, we have 
%\begin{align*}\widehat{\Omega}=U^{*}\Omega=\begin{pmatrix}
%U_{1}^{*}\Omega \\ 
%U_{2}^{*}\Omega
%\end{pmatrix} =\begin{pmatrix}
%\widehat{\Omega}_{1}\\ 
%\widehat{\Omega}_{2}
%\end{pmatrix} ,
%\end{align*}
%where $ \widehat{\Omega}_{1}=U_{1}^{*}\Omega\in\mathbb{C}^{l-p\times l}$,  $\widehat{\Omega}_{2}=U_{2}^{*}\Omega\in\mathbb{C}^{n-l+p\times l}$, $0\leq p\leq l-k$, parameter $p$ is to balance the need for oversampling for reliability.\\
%We assume $ {\rm rank}\left(\widehat{\Omega}_{1}\right)=l-p $, then its pseudo-inverse satisfies \begin{align*}\widehat{\Omega}_{1}\widehat{\Omega}_{1}^{\dagger }=I_{l-p}.\end{align*}
Using the notation introduced in Section \ref{sec:2.1}, we write $K$ as follows,
\begin{align}\label{3000}
K= U\phi\left(\Sigma \right )\widehat{\Omega}=U\begin{pmatrix}
\phi\left( \Lambda_{1}\right) &  & \\ 
& \phi\left( \Lambda _{2}\right)
\end{pmatrix}\begin{pmatrix}
\widehat{\Omega}_{1}\\ 
\widehat{\Omega}_{2}
\end{pmatrix}=U\begin{pmatrix}
\phi\left(\Lambda_{1}\right)\widehat{\Omega}_{1}\\
\phi\left(\Lambda_{2}\right)\widehat{\Omega}_{2}
\end{pmatrix}.
\end{align}
To make the last block matrix in \eqref{3000} be simplified, we choose a matrix $ X $ in the following form:
\begin{align*}
X= \begin{pmatrix}
\widehat{\Omega}_{1}^{\dagger}\phi^{-1}\left( \Lambda _{1}\right), & X_{2}
\end{pmatrix}\in\mathbb{C}^{l\times l},
\end{align*}
where $ X_{2}\in \mathbb{C}^{l\times p} $ is such that $ X $ is nonsingular and $ \widehat{\Omega}_{1}X_{2}=0 $. %Now we find the X satisfy the conditon of Lemma \ref{lem201}.\\
%This constructed $ X $ multiplied by $K$ leads to 
Thus, 
\begin{align*}
KX=U\begin{pmatrix}
I_k& 0  \\ 
H_{1}& H_{2} 
\end{pmatrix},
\end{align*}
where
\begin{align}\label{302}
H_{1}=\phi\left(\Lambda _{2}\right)\widehat{\Omega}_{2}\widehat{\Omega}_{1}^{\dagger}\phi^{{-1}}\left(\Lambda _{1}\right),\quad H_{2}=\phi\left(\Lambda _{2}\right)\widehat{\Omega}_{2}X_{2}.
\end{align}
In accord with the above block form, we write the thin $QR$ factorization of $KX$ in the following form:
\begin{align}\label{303}
KX=\widetilde{Q}\widetilde{R}=\begin{pmatrix}
\widetilde{Q}_{1}, & \widetilde{Q}_{2}
\end{pmatrix}\begin{pmatrix}
\widetilde{R}_{11}& \widetilde{R}_{12}  \\ 
&   \widetilde{R}_{22}
\end{pmatrix}=U\begin{pmatrix}
I_k& 0 \\ 
H_{1}& H_{2} 
\end{pmatrix}.
\end{align}
As a result, we have the following thin QR factorization %Actually, there is a small $QR$ factorization in \eqref{303} which we will use in the following analysis,
\begin{align}\label{304}
U\begin{pmatrix}
I_k\\ 
H_{1}
\end{pmatrix}=\widetilde{Q}_{1}\widetilde{R}_{11},
\end{align}
which will be used in our error analysis.
%Below we present deterministic analysis for the trace estimator and log-determiant estimator.

%To prove the error bounds for log-determinant, 
Furthermore, we also need the following three results.
\begin{lemma}\label{lem203} %{\rm\cite{halko}}
	Suppose ${\rm range}(N) \subset {\rm range}(M)$. Then, for any Hermitian positive semi-definite matrix $A$, the following inequality holds
	\begin{align*}
	{\rm Tr}(P_NA)\leqslant{\rm Tr}(P_MA),
	\end{align*}
	where $P_N$ and $P_M$ are the orthogonal projections on ${\rm range}(N)$ and ${\rm range}(M)$, respectively.
\end{lemma}
\emph{Proof} From the proof of \cite[Proposition 8.5]{halko}, we know $P_N\preceq P_M$, where $\preceq$ denotes the L\"owner partial order \cite[Definition 7.7.1]{horn13}. Then, by the known conjugation rule (see e.g., \cite[Theorem 7.7.2]{horn13}),
$$A^{\frac{1}{2}}P_NA^{\frac{1}{2}}\preceq A^{\frac{1}{2}}P_MA^{\frac{1}{2}}. $$
Further, by the properties of L\"owner partial order and trace, we have
$${\rm Tr}\left(A^{\frac{1}{2}}P_NA^{\frac{1}{2}}\right)\preceq {\rm Tr}\left(A^{\frac{1}{2}}P_MA^{\frac{1}{2}}\right),\ \textrm{ i.e.,}\  {\rm Tr}(P_NA)\leqslant{\rm Tr}(P_MA). $$ \qed
\begin{lemma}\label{lem2030} 
	For any two Hermitian positive semi-definite matrices $A$ and $B$ with the same order, the following inequalities hold
	\begin{align*}
	{\rm Tr}(AB)\leqslant{\rm Tr}(A)\lambda_{\max}(B)\leqslant{\rm Tr}(A){\rm Tr}(B),\\
	\quad{\rm Tr}(AB)\leqslant\lambda_{\max}(A){\rm Tr}(B)\leqslant{\rm Tr}(A){\rm Tr}(B),
	\end{align*}
	where $\lambda_{\max}(A)$ and $\lambda_{\max}(B)$ are the largest eigenvalues of $A$ and $B$, respectively.
\end{lemma}
\emph{Proof} These inequalities are well-known results and can be derived from, e.g., von Neumann trace theorem \cite[Theorem 7.4.1.1]{horn13} directly. \qed

\begin{lemma}\label{lem202}{\rm\cite[Corollary 2.1]{ouell}}
	If $ B\in \mathbb{C}^{m \times n} $ and $ C\in \mathbb{C}^{n \times m} $, then
	\begin{align*}
	\det(I_{m}\pm BC)=\det(I_{n}\pm CB).
	\end{align*}
\end{lemma}

\subsubsection{Structural bounds for trace estimator}
\label{sec:3.3.1}
%We first prove the following structural (deterministic) bounds, i.e., we temporarily consider $\Omega$ to be a deterministic matrix. %Then, we derive probabilistic bounds 
\begin{theorem}\label{them310}
	Let $T=Q_{q}^{*}AQ_{q}$ be computed by Algorithm \ref{Al2}. Then
	\begin{align}\label{10000}
	0\leq {\rm Tr}(A)-{\rm Tr}(T)\leq \left(1+T_{q-1}^{-1}\left(\frac{2\lambda_{k}-\lambda_{k+1}}{\lambda_{k+1}}\right)\left \| \widehat{\Omega}_{2}\widehat{\Omega}_{1}^{\dagger}\right\|_{2}\right) {\rm Tr}\left( \Lambda _{2}\right ).
	\end{align}
	When $0<\left \| \widehat{\Omega}_{2}\widehat{\Omega}_{1}^{\dagger}\right\|_{2}\leqslant\frac{\lambda _{k}}{\lambda _{k+1}}T_{q-1}(\frac{2\lambda_{k}-\lambda_{k+1}}{\lambda_{k+1}})$, the following bound is tighter, % with $ \gamma =\frac{2\lambda_{k}-\lambda_{k+1}}{\lambda_{k+1}}$. %\geq \frac{\lambda_{k}}{\lambda_{l-p+1}}. $
	\begin{align}\label{10001}
	0\leq {\rm Tr}(A)-{\rm Tr}(T)\leq \left(1+\frac{\lambda _{k+1}}{\lambda _{k}}T_{q-1}^{-2}\left(\frac{2\lambda_{k}-\lambda_{k+1}}{\lambda_{k+1}}\right)\left \| \widehat{\Omega}_{2}\widehat{\Omega}_{1}^{\dagger}\right\|_{2}^{2}\right) {\rm Tr}\left( \Lambda _{2}\right ).
	\end{align}
\end{theorem}
\emph{Proof} The lower bound has been proven in \cite[Lemma 1]{saibaba}. In the following, we show that the upper bounds hold. % is proved as follows.\\

Since $K= \phi\left(A\right)\Omega$ is an element of $\mathcal{K}_q$, we have 
$${\rm range}(K)\subset  \mathcal{K}_q.$$
Thus, by Lemma \ref{lem201}, we get 
\begin{align}\label{305}
{\rm range}(\widetilde{Q})={\rm range}(Q)\subseteq {\rm range}(Q_{q}),
\end{align}
which together with Lemma \ref{lem203} and \eqref{303} implies
\begin{eqnarray}
\label{306}	
{\rm Tr}(A)-{\rm Tr}(T)&=&{\rm Tr}(A)-{\rm Tr}(Q_{q}^{*}AQ_{q})={\rm Tr}(A)-{\rm Tr}(Q_{q}Q_{q}^{*}A)\nonumber\\
&\leq& {\rm Tr}(A)-{\rm Tr}(QQ^{*}A)\nonumber
={\rm Tr}(A)-{\rm Tr}(\widetilde{Q}\widetilde{Q}^{*}A)\\
&\leq& {\rm Tr}(A)-{\rm Tr}(\widetilde{Q}_{1}\widetilde{Q}_{1}^{*}A)={\rm Tr}(A)-{\rm Tr}(\widetilde{Q}_{1}^{*}A\widetilde{Q}_{1}).
\end{eqnarray}
From \eqref{304}, we have 
\begin{align}\label{307}
\widetilde{Q}_{1}=U\begin{pmatrix}
I_k\\ 
H_{1}
\end{pmatrix}\widetilde{R}_{11}^{-1}\ \textrm{  and  }\ \widetilde{R}_{11}^{*}\widetilde{R}_{11}=I_k+H_{1}^{*}H_{1}.
\end{align}
As a result, 
\begin{eqnarray}\label{308}
\widetilde{Q}_{1}^{*}A\widetilde{Q}_{1}&=&(\widetilde{R}_{11}^{*})^{-1}\begin{pmatrix}
I_k,  &H_{1}^{*} 
\end{pmatrix}U^{*}U\begin{pmatrix}
\Lambda_{1} &  \\ 
&\Lambda_{2} 
\end{pmatrix}U^{*}U\begin{pmatrix}
I_k\\ 
H_{1}
\end{pmatrix}\widetilde{R}_{11}^{-1}\nonumber\\
&=&(\widetilde{R}_{11}^{*})^{-1}(\Lambda_{1}+H_{1}^{*}\Lambda_{2}H_{1})\widetilde{R}_{11}^{-1},
\end{eqnarray}
and hence
\begin{eqnarray} \label{309}
{\rm Tr}(\widetilde{Q}_{1}^{*}A\widetilde{Q}_{1})&=&{\rm Tr}((\Lambda_{1}+H_{1}^{*}\Lambda_{2}H_{1})(\widetilde{R}_{11}^{*}\widetilde{R}_{11})^{-1})
={\rm Tr}((\Lambda_{1}+H_{1}^{*}\Lambda_{2}H_{1})(I_k+H_{1}^{*}H_{1})^{-1})\nonumber\\
&=&{\rm Tr}(\Lambda_{1}(I_k+H_{1}^{*}H_{1})^{-1})+{\rm Tr}(\Lambda_{2}H_{1}(I_k+H_{1}^{*}H_{1})^{-1}H_{1}^{*}).
\end{eqnarray}
Substituting \eqref{309} into \eqref{306} and noting ${\rm Tr}(A)={\rm Tr}(\Lambda_{1})+{\rm Tr}(\Lambda_{2})$ and Lemma \ref{lem2030} gives
\begin{eqnarray}\label{310} 
{\rm Tr}(A)-{\rm Tr}(T)&\leq&{\rm Tr}(\Lambda_{1}(I_k-(I_k+H_{1}^{*}H_{1})^{-1}))+{\rm Tr}(\Lambda_{2}(I_k-H_{1}(I_k+H_{1}^{*}H_{1})^{-1}H_{1}^{*}))\nonumber\\
&\leq& {\rm Tr}(\Lambda_{1}(I_k-(I_k+H_{1}^{*}H_{1})^{-1}))\nonumber\\
&+&{\rm Tr}\begin{pmatrix}
\Lambda_{2}
\end{pmatrix}\lambda_{\max}(I_k-H_{1}(I_k+H_{1}^{*}H_{1})^{-1}H_{1}^{*})\nonumber\\
&\leq& {\rm Tr}(\Lambda_{1}(I_k-(I_k+H_{1}^{*}H_{1})^{-1}))+{\rm Tr}\begin{pmatrix}
\Lambda_{2}
\end{pmatrix}.
\end{eqnarray}
Note that 
\begin{align*}
(I_{k}+H_{1}^{*}H_{1})^{-1}=I_{k}-H_{1}^{*}(I_{n-k}+H_{1}H_{1}^{*})^{-1}H_{1}.
\end{align*}
Then 
$${\rm Tr}(\Lambda_{1}(I_{k}-(I_{k}+H_{1}^{*}H_{1})^{-1}))={\rm Tr}(\Lambda_{1}H_{1}^{*}(I_{n-k}+H_{1}H_{1}^{*})^{-1}H_{1}).$$
Further, setting $\phi(x)=xf(x)$, where $f(x)$ is defined in \eqref{202}, and considering $H_{1}$ in \eqref{302}, von Neumann trace theorem \cite[Theorem 7.4.1.1]{horn13}, and singular value inequalities \cite[Theorem 3.3.14]{horn}, we have %the first term in the last inequality of \eqref{310} can have the following form,
\begin{eqnarray*}
	{\rm Tr}(\Lambda_{1}(I_{k}&-&(I_{k}+H_{1}^{*}H_{1})^{-1}))%&=&{\rm Tr}(\Lambda_{1}\phi^{-1}(\Lambda _{1})(\widehat{\Omega}_{2}\widehat{\Omega}_{1}^{\dagger})^{*}\phi(\Lambda _{3})(I_{n-k}+H_{1}H_{1}^{*})^{-1}H_{1})\nonumber\\
	\\&=&{\rm Tr}(\Lambda_{1}\Lambda_{1}^{-1}f^{-1}(\Lambda _{1})(\widehat{\Omega}_{2}\widehat{\Omega}_{1}^{\dagger})^{*}\Lambda_{2}f(\Lambda _{2})(I_{n-k}+H_{1}H_{1}^{*})^{-1}H_{1})\nonumber\\
	&\leq& \sum_{j=1}^{k}\sigma_{j}\left(f^{-1}(\Lambda _{1})(\widehat{\Omega}_{2}\widehat{\Omega}_{1}^{\dagger})^{*}\Lambda_{2}f(\Lambda _{2})\right)\sigma_{j}\left( (I_{n-k}+H_{1}H_{1}^{*})^{-1}H_{1}\right )\nonumber\\
	&\leq& \sum_{j=1}^{k}\sigma_{j}\left(f^{-1}(\Lambda _{1})(\widehat{\Omega}_{2}\widehat{\Omega}_{1}^{\dagger})^{*}\Lambda_{2}f(\Lambda _{2})\right)\left \| (I_{n-k}+H_{1}H_{1}^{*})^{-1}H_{1} \right \|_{2}\nonumber\\
	&\leq& \left \| f^{-1}(\Lambda _{1}) \right \|_{2}\left \| f(\Lambda_{2})\right\|_{2}\left \| \widehat{\Omega}_{2}\widehat{\Omega}_{1}^{\dagger} \right \|_{2}\left \| (I_{n-k}+H_{1}H_{1}^{*})^{-1}H_{1} \right \|_{2}{\rm Tr}(\Lambda_{2}),
\end{eqnarray*}
which combined with \eqref{203} and \eqref{204} leads to
\begin{eqnarray}\label{311} 
&&{\rm Tr}(\Lambda_{1}(I_{k}-(I_{k}+H_{1}^{*}H_{1})^{-1}))
\nonumber\\
&&\leq T_{q-1}^{-1}\left(\frac{2\lambda _{k}-\lambda _{k+1}}{\lambda _{k+1}}\right)\left \| \widehat{\Omega}_{2}\widehat{\Omega}_{1}^{\dagger} \right \|_{2}\left \| (I_{n-k}+H_{1}H_{1}^{*})^{-1}H_{1} \right \|_{2}{\rm Tr}(\Lambda_{2}).
\end{eqnarray}
Furthermore, from \cite{saibaba}, we have
\begin{align}\label{312}
\left \| (I_{n-k}+H_{1}H_{1}^{*})^{-1}H_{1} \right \|_{2}\leq 1, 
\end{align}
or  
\begin{eqnarray}\label{313}
\left \| (I_{n-k}+H_{1}H_{1}^{*})^{-1}H_{1} \right \|_{2}&\leq& \left\|H_{1}\right\|_{2}\leq \left\|\Lambda _{2}f(\Lambda_{2}) \right\|_{2}\left\|\Lambda_{1}^{-1}f^{-1}(\Lambda_{1}) \right\|_{2}\left\|\widehat{\Omega}_{2}\widehat{\Omega}_{1}^{\dagger}\right\|_{2}\nonumber \\
&\leq& \frac{\lambda_{k+1}}{\lambda_{k}}T_{q-1}^{-1}\left(\frac{2\lambda _{k}-\lambda _{k+1}}{\lambda _{k+1}}\right)\left\| \widehat{\Omega}_{2}\widehat{\Omega}_{1}^{\dagger}\right\|_{2}.
\end{eqnarray}
Thus, combining \eqref{310}, \eqref{311}, \eqref{312}, and \eqref{313}, we derive the desired upper bounds \eqref{10000} and \eqref{10001}. \qed
%Using the following result \cite[Lemma 4]{saibaba} and as done in \cite[Section 4.1.1]{saibaba}, we can complete the proof of Theorem \ref{lem308}.
\subsubsection{Structural bounds for log-determiant estimator}\label{Sec3.3.2}
\begin{theorem}\label{them311}
	Let $T=Q_{q}^{*}AQ_{q}$ be computed by Algorithm \ref{Al2}. Then
	\begin{eqnarray}\label{40000}
	0&\leq& \log\det(I_{n}+A)-\log\det(I_{ql}+T)\nonumber\\
	&\leq& \log\det(I_{n-k}+\eta\Lambda_{2})+\log\det(I_{n-k}+
	\Lambda_{2} ),
	\end{eqnarray}
	where $\eta=\frac{\lambda _{k+1}}{\lambda _{k}}T_{q-1}^{-2}(\frac{2\lambda_{k}-\lambda_{k+1}}{\lambda_{k+1}})\left \| \widehat{\Omega}_{2}\widehat{\Omega}_{1}^{\dagger}\right\|_{2}^{2}$. %, and $ \gamma =\frac{2\lambda_{l-p}-\lambda_{l-p+1}}{\lambda_{l-p+1}}\geq \frac{\lambda_{l-p}}{\lambda_{l-p+1}}. $
\end{theorem}
\emph{Proof}   The lower bound has been derived in \cite[Lemma 2]{saibaba}. It is sufficient to show that the upper bound holds.
%Similar to formula \eqref{306}, since \eqref{305} is established, we can derive
%By Lemma \ref{lem202} and \eqref{305}, we have

From \eqref{305}, it follows that there is an orthonormal matrix $Y\in \mathbb{C}^{l \times l}$ such that
$Q=Q_qY$, and hence
$$
Q^*AQ=Y^*Q_q^*AQ_qY=Y^*TY.
$$
Thus, by the proved lower bound, i.e., 
\begin{align}\label{30000}
\log\det(I_{n}+A)-\log\det(I_{ql}+T)\geqslant 0,
\end{align}
we have
\begin{align*}
\log\det(I_{ql}+T)-\log\det(I_{l}+Y^*TY)\geqslant 0.
\end{align*}
That is,
\begin{align*}
\log\det(I_{ql}+Q_q^*AQ_q)-\log\det(I_{l}+Q^*AQ)\geqslant 0.
\end{align*}
Thus, 
\begin{equation*}
\log\det(I_{n}+A)-\log\det(I_{ql}+T)\leq \log\det(I_{n}+A)-\log\det(I_{l}+Q^{*}AQ).
\end{equation*}
By Lemmas \ref{lem202} and \ref{lem201}, it is seen that
\begin{eqnarray}\label{314}
\log\det(I_{n}+A)-\log\det(I_{ql}+T)&\leq& \log\det(I_{n}+A)-\log\det(I_{l}+QQ^{*}A)\nonumber\\
&=&\log\det(I_{n}+A)-\log\det(I_{l}+\widetilde{Q}\widetilde{Q}^{*}A)\nonumber\\
&=&\log\det(I_{n}+A)-\log\det(I_{l}+\widetilde{Q}^{*}A\widetilde{Q}).
\end{eqnarray}
From \eqref{303}, we can check that 
$\widetilde{Q}_{1}=\widetilde{Q}\begin{pmatrix}
I_k\\ 
0
\end{pmatrix}$
and $\begin{pmatrix}
I_k\\ 
0
\end{pmatrix}$ is orthonormal.  Thus, by \eqref{30000} again, we have
\begin{eqnarray*}
	\log\det(I_{l}+\widetilde{Q}^{*}A\widetilde{Q})\geqslant\log\det(I_{k}+\widetilde{Q}_1^{*}A\widetilde{Q}_1),
\end{eqnarray*}
which together with \eqref{314} gives 
\begin{align}\label{3150}
\log\det(I_{n}+A)-\log\det(I_{ql}+T)&\leq\log\det(I_{n}+A)-\log\det(I_{k}+\widetilde{Q}_{1}^{*}A\widetilde{Q}_{1}).
\end{align}
%The reasoning process above has used Lemma \ref{lem202}, the places where Lemma \ref{lem202} is used below, no longer mentioned.\\
Setting $M=\widetilde{Q}_{1}^{*}A\widetilde{Q}_{1}$ and considering \eqref{308}, %\eqref{201}, and \eqref{20000}, 
we get
\begin{align*}
M=(\widetilde{R}_{11}^{*})^{-1}(\Lambda_{1}+H_{1}^{*}\Lambda_{2}H_{1})\widetilde{R}_{11}^{-1}.
\end{align*}
%And define function $ h(\cdot) $ as
%\begin{align*}
%h(\cdot)=\log\det(\cdot).
%\end{align*}
%Then we have
Thus, by Lemma \ref{lem202} and noting \eqref{307}, we have
\begin{align}\label{30001}
\log\det(I_{k}+M)=\log\det(I_{k}+M_{1})=\log\det(I_{k}+M_{2}),
\end{align}
where 
\begin{align*}
M_{1}&=(\Lambda_{1}+H_{1}^{*}\Lambda_{2}H_{1})(\widetilde{R}_{11}^{*}\widetilde{R}_{11})^{-1}=(\Lambda_{1}+H_{1}^{*}\Lambda_{2}H_{1})(I+H_{1}^{*}H_{1})^{-1},\\
M_{2}&=(I+H_{1}^{*}H_{1})^{-\frac{1}{2}}(\Lambda_{1}+H_{1}^{*}\Lambda_{2}H_{1})(I+H_{1}^{*}H_{1})^{-\frac{1}{2}}\\
&\succeq (I+H_{1}^{*}H_{1})^{-\frac{1}{2}}\Lambda_{1}(I+H_{1}^{*}H_{1})^{-\frac{1}{2}}=M_{3}.
\end{align*}
Combining the properties of L\"owner partial order \cite[Corollary 7.7.4]{horn13} with \eqref{30001} implies
%\begin{align*}
%I_{k}+M_{2}\geq I_{k}+M_{3},
%\end{align*}
%which lead to 
\begin{align*}
\log\det(I_{k}+M)\geq \log\det(I_{k}+M_{3}).
\end{align*}
As done in \cite{saibaba},  we can show that $\log\det(I_{k}+M_{3})= \log\det(I_{k}+M_{4})$, where
\begin{align*}
M_{4}{=}\Lambda_{1}^{\frac{1}{2}}(I+H_{1}^{*}H_{1})^{-1}\Lambda_{1}^{\frac{1}{2}},
\end{align*}
and 
\begin{align*}
\log\det(I_{k}+\Lambda_{1})-\log\det(I_{k}+M_{4})=\log\det(M_{5}),
\end{align*}
where 
\begin{align*}
M_{5}&=(I_{k}+M_{4})^{-\frac{1}{2}}(I_{k}+\Lambda_{1})(I_{k}+M_{4})^{-\frac{1}{2}}\\
&=(I_{k}+M_{4})^{-1}+(I_{k}+M_{4})^{-\frac{1}{2}}\Lambda_{1}(I_{k}+M_{4})^{-\frac{1}{2}}.
\end{align*}
Further, as done in \cite{saibaba}, we have
\begin{align*}
M_{5}=I_{k}+(I_{k}+M_{4})^{-\frac{1}{2}}(\Lambda_{1}-M_{4})(I_{k}+M_{4})^{-\frac{1}{2}}=I_{k}+M_{6},
\end{align*}
where
\begin{align*}
M_{6}=(I_{k}+M_{4})^{-\frac{1}{2}}(\Lambda_{1}-M_{4})(I_{k}+M_{4})^{-\frac{1}{2}}\preceq(I_{k}+M_{4})^{-\frac{1}{2}}\Lambda_{1}^{\frac{1}{2}}H_{1}^{*}H_{1}\Lambda_{1}^{\frac{1}{2}}(I_{k}+M_{4})^{-\frac{1}{2}},
\end{align*}
and hence
%From \eqref{316} we have 
\begin{align*}
&M_{5}\preceq I_{k}+(I_{k}+M_{4})^{-\frac{1}{2}}\Lambda_{1}^{\frac{1}{2}}H_{1}^{*}H_{1}\Lambda_{1}^{\frac{1}{2}}(I_{k}+M_{4})^{-\frac{1}{2}}\preceq I_{k}+\Lambda_{1}^{\frac{1}{2}}H_{1}^{*}H_{1}\Lambda_{1}^{\frac{1}{2}}.
\end{align*}
Thus, considering \eqref{3150} and \eqref{30001} and using Lemma \ref{lem202}, we have
\begin{eqnarray}\label{317}
\log\det(I_{n}+A)&-&\log\det(I_{ql}+T)\nonumber\\
&\leq& \log\det(I_{k}+\Lambda_{1} )-\log\det(I_{k}+M)+\log\det(I_{n-k}+\Lambda_{2} )\nonumber\\
&\leq&\log\det(I_{k}+\Lambda_{1}^{\frac{1}{2}}H_{1}^{*}H_{1}\Lambda_{1}^{\frac{1}{2}})+\log\det(I_{n-k}+\Lambda_{2} )\nonumber\\
&=&\log\det(I_{n-k}+H_{1}\Lambda_{1}H_{1}^{*})+\log\det(I_{n-k}+
\Lambda_{2} ).
\end{eqnarray}
Further, noting $H_1$ in \eqref{302}, $\phi(x)=xf(x)$, and Lemma \ref{lem202}, we have
\begin{eqnarray}\label{318}
\quad  &&\log\det(I_{n-k}+H_{1}\Lambda_{1}H_{1}^{*})\nonumber\\
&&=\log\det\left(I_{n-k}+\phi(\Lambda _{2})\widehat{\Omega}_{2}\widehat{\Omega}_{1}^{\dagger}
\phi^{-1}(\Lambda_{1})\Lambda_{1}
\phi^{-1}(\Lambda_{1})(\widehat{\Omega}_{2}\widehat{\Omega}_{1}^{\dagger})^{*}\phi(\Lambda_{2})\right)\nonumber\\
&&=\log\det(I_{n-k}+\Lambda_{2}f(\Lambda_{2})\widehat{\Omega}_{2}\widehat{\Omega}_{1}^{\dagger}
f^{-1}(\Lambda_{1})\Lambda_{1}^{-1}
f^{-1}(\Lambda_{1})(\widehat{\Omega}_{2}\widehat{\Omega}_{1}^{\dagger})^{*}\Lambda_{2}f(\Lambda_{2}))\nonumber\\
&&=\log\det(I_{k}+\Lambda_{1}^{-\frac{1}{2}}
f^{-1}(\Lambda_{1})(\widehat{\Omega}_{2}\widehat{\Omega}_{1}^{\dagger})^{*}f(\Lambda_{2})\Lambda_{2}^{\frac{1}{2}}\Lambda_{2}\Lambda_{2}^{\frac{1}{2}}f(\Lambda_{2})\widehat{\Omega}_{2}\widehat{\Omega}_{1}^{\dagger}
f^{-1}(\Lambda_{1})\Lambda_{1}^{-\frac{1}{2}})\nonumber\\
&&=\log\det(I_{k}+G^{*}\Lambda_{2}G),
\end{eqnarray}
where 
\begin{align*}
G=\Lambda_{2}^{\frac{1}{2}}f(\Lambda_{2})\widehat{\Omega}_{2}\widehat{\Omega}_{1}^{\dagger}
f^{-1}(\Lambda_{1})\Lambda_{1}^{-\frac{1}{2}}.
\end{align*}
As done in \cite{saibaba}, from \cite[Theorem 3.3.16]{horn}, we can derive 
\begin{equation*}%\%label{319}
\begin{split}
\det(I_{k}+G^{*}\Lambda_{2}G)&\leq \prod_{j=1}^{n-k}(I_{n-k}+\sigma_{j}(GG^{*}\Lambda_{2}))\leq \prod_{j=1}^{n-k}(I_{n-k}+\left\| GG^{*}\right\|_{2}\sigma_{j}(\Lambda_{2}))\\
&\leq \det(I_{n-k}+\left\| G\right\|_{2}^{2}\Lambda_{2}).
\end{split}
\end{equation*}
Substituting the about result into \eqref{318} and then into \eqref{317} gives
\begin{eqnarray}\label{319}
\log\det(I_{n}+A)&-&\log\det(I_{ql}+T)\nonumber\\
&\leq& \log\det(I_{n-k}+\left\| G\right\|_{2}^{2}\Lambda_{2})+\log\det(I_{n-k}+
\Lambda_{2} ).
\end{eqnarray}
Further, noting \eqref{203} and \eqref{204}, we obtain
\begin{equation*}%\label{320}
\begin{split}
\left\| G\right\|_{2}^{2}=\left\| \Lambda_{2}^{\frac{1}{2}}f(\Lambda_{2})\widehat{\Omega}_{2}\widehat{\Omega}_{1}^{\dagger}
f^{-1}(\Lambda _{1})\Lambda_{1}^{-\frac{1}{2}}\right\|_{2}^{2}&\leq \frac{\lambda_{k+1}}{\lambda_{k}}\left \| f^{-1}(\Lambda_{1}) \right \|_{2}^{2}\left \| f(\Lambda_{2})\right\|_{2}^{2}\left \| \widehat{\Omega}_{2}\widehat{\Omega}_{1}^{\dagger}\right \|_{2}^{2}\\
&\leq \frac{\lambda_{k+1}}{\lambda_{k}}T_{q-1}^{-2}\left(\frac{2\lambda _{k}-\lambda _{k+1}}{\lambda _{k+1}}\right)\left\| \widehat{\Omega}_{2}\widehat{\Omega}_{1}^{\dagger}\right\|_{2}^{2}\\&=\eta,
\end{split}
\end{equation*}
which together with \eqref{319} implies the desired upper bound.  \qed %Then we can combine \eqref{317}, \eqref{318}, \eqref{319}, \eqref{320}, come to the result
%\begin{align*}
%\log\det(I_{n}+A)-\log\det(I_{ql}+T)&\leq \log\det(I_{n-l+p}+\left\| G\right\|_{2}^{2}\Lambda_{3})+\log\det(I_{n-k}+\begin{pmatrix}
%\Lambda_{2} & \\ 
%& \Lambda_{3}
%\end{pmatrix})\\
%&\leq \log\det(I_{n-l+p}+\eta_{2}\Lambda_{3})+\log\det(I_{n-k}+\begin{pmatrix}
%\Lambda_{2} & \\ 
%& \Lambda_{3}
%\end{pmatrix}).
%\end{align*}

%\subsubsection{Probabilistic error bounds }\label{Sec3.3.3}

\section{Numerical experiments}
\label{sec:4}
In this section, we take two examples from \cite{saibaba} to demonstrate the performance of our algorithms and compare them with randomized subspace iteration algorithms given in \cite{saibaba}. We also test the structural upper bounds given Theorems \ref{them310} and \ref{them311} by these two examples.  In these experiments, the relative errors in trace and log-determinant estimators are defined as
\begin{align*}
\Delta_{t}\equiv\frac{{\rm Tr}(A)-{\rm Tr}(T)}{{\rm Tr}(A)},\quad
\Delta_{l}\equiv\frac{\log\det(I+A)-\log\det(I+T)}{\log\det(I+A)},
\end{align*}%For convenience, our experiments are based on two-bolck patition of $ \Sigma $. 
and all computations are carried out in MATLAB 2016b.

\subsection{Small matrices}\label{Sec5.1}

%In this section we consider the performance of our algorithm on small test examples.
The eigenvalues of the test matrix $ A $ satisfy $ \lambda_{j+1}=\tau^{j}\lambda_{1} $ for $ j=1,2,\cdots,n-1$. In contrast to \cite{saibaba}, we set the order of the matrix $A$ to be $ 1280 \times 1280$. %, which is quite small in big data age. 
By setting suitable values of $\tau$ and $\lambda_{1}$, % and different  $\lambda_{1}$. , 
we do the following four specific numerical experiments.
\begin{enumerate}
	\item Test the performance of Algorithm \ref{Al2} when $p = 20$, $q = 3$, $\lambda_{1} = 100$, and $\tau$ varies from 0.98 to 0.86.
	\item Test the performance of Algorithm \ref{Al2} when $p = 20$, $\lambda_{1} = 100$, $\tau = 0.90$, and $q$ varies from 1 to 5.
	\item Compare Algorithms \ref{Al1} and \ref{Al2}, when $p = 20$, $q = 3$, $\lambda_{1} = 100$, and $\tau = 0.92$.
	\item Test the structural error bounds when  $p = 20$, $q = 3$, $\lambda_{1} = 100$, and $\tau = 0.90$. %, and compare the
\end{enumerate}

\begin{figure}[H]
	\centering
	\subfigure{%[SubfigureCaption]{
		%\label{Fig.sub.1}
		\includegraphics[height=48mm,width=55mm]{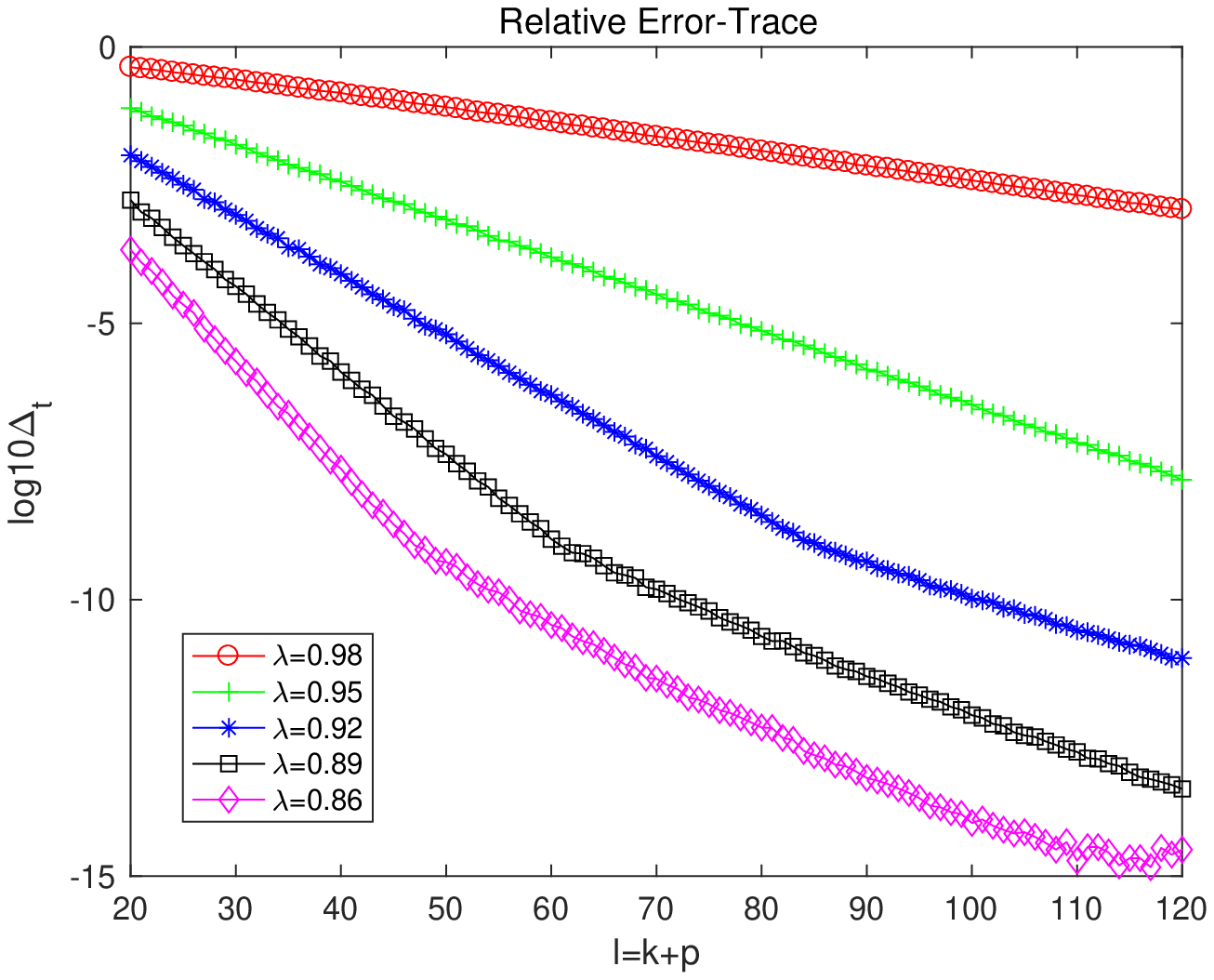}}
	\subfigure{%[SubfigureCaption]{
		%\label{Fig.sub.2}
		\includegraphics[height=48mm,width=55mm]{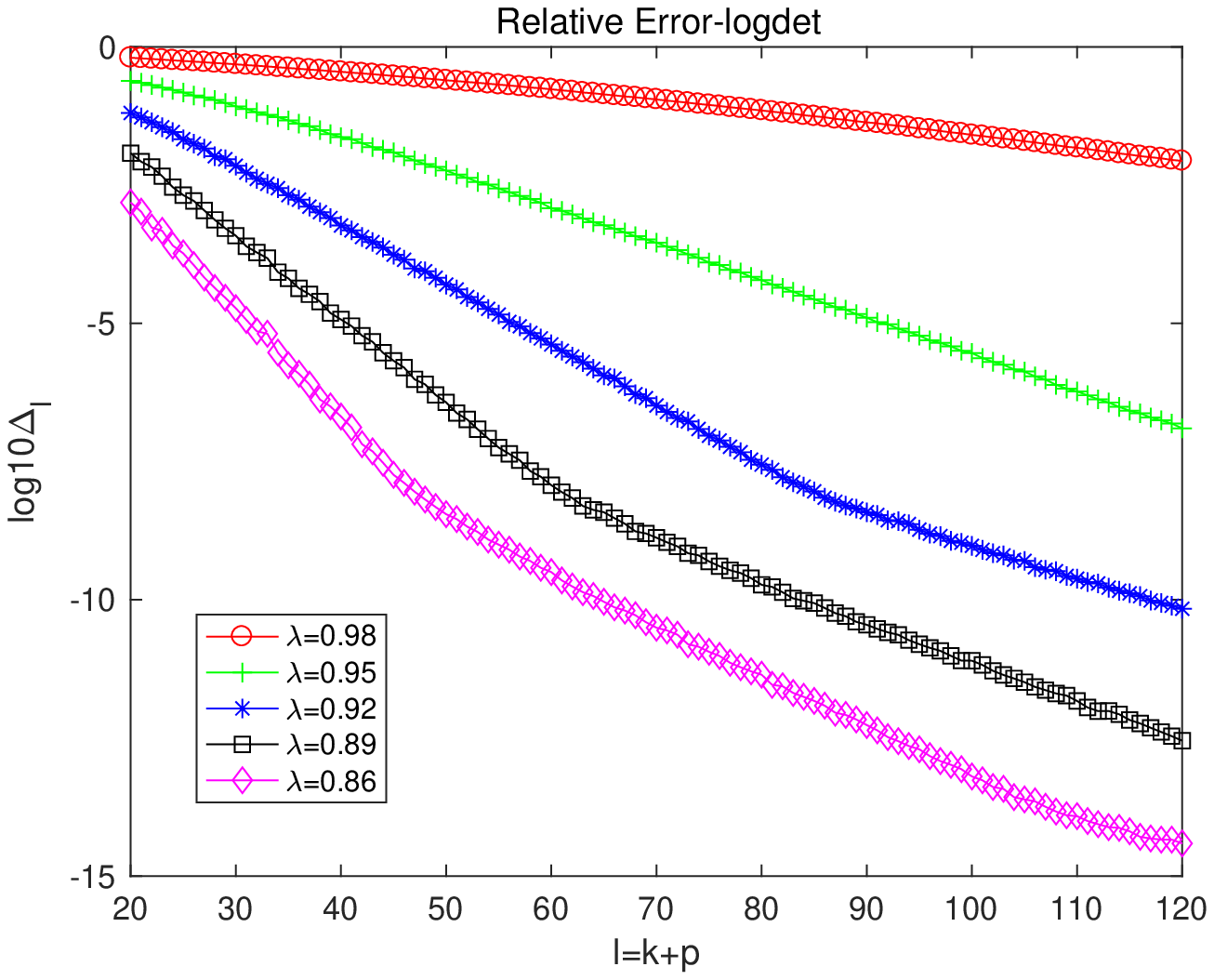}}
	\caption{%Accuracy of proposed estimators on a matrix with geometrically decaying eigenvalues.
	Accuracy of (left) trace and (right) log-determinant estimators produced by Algorithm \ref{Al2} for small matrix with $\tau$ varying from $0.98$ to $0.86$. The relative error is plotted against the sample size}
	\label{Fig1}
\end{figure}

The first experiment is used to test the effect of gap on the algorithms. Numerical results are displayed in Fig. \ref{Fig1}. It is easy to see that both the trace and log-determinant estimators are increasingly accurate as the eigenvalue gap increases. Note that hereafter the relative error is plotted against the sample size $l=k+p$. Since $p$ is fixed, increasing the sample size means to increase in the target rank $k$. As a result, the location of the gap is changing. 

\begin{figure}[H]
	\centering
	\subfigure{%[SubfigureCaption]{
		%\label{Fig.sub.1}
		\includegraphics[height=48mm,width=55mm]{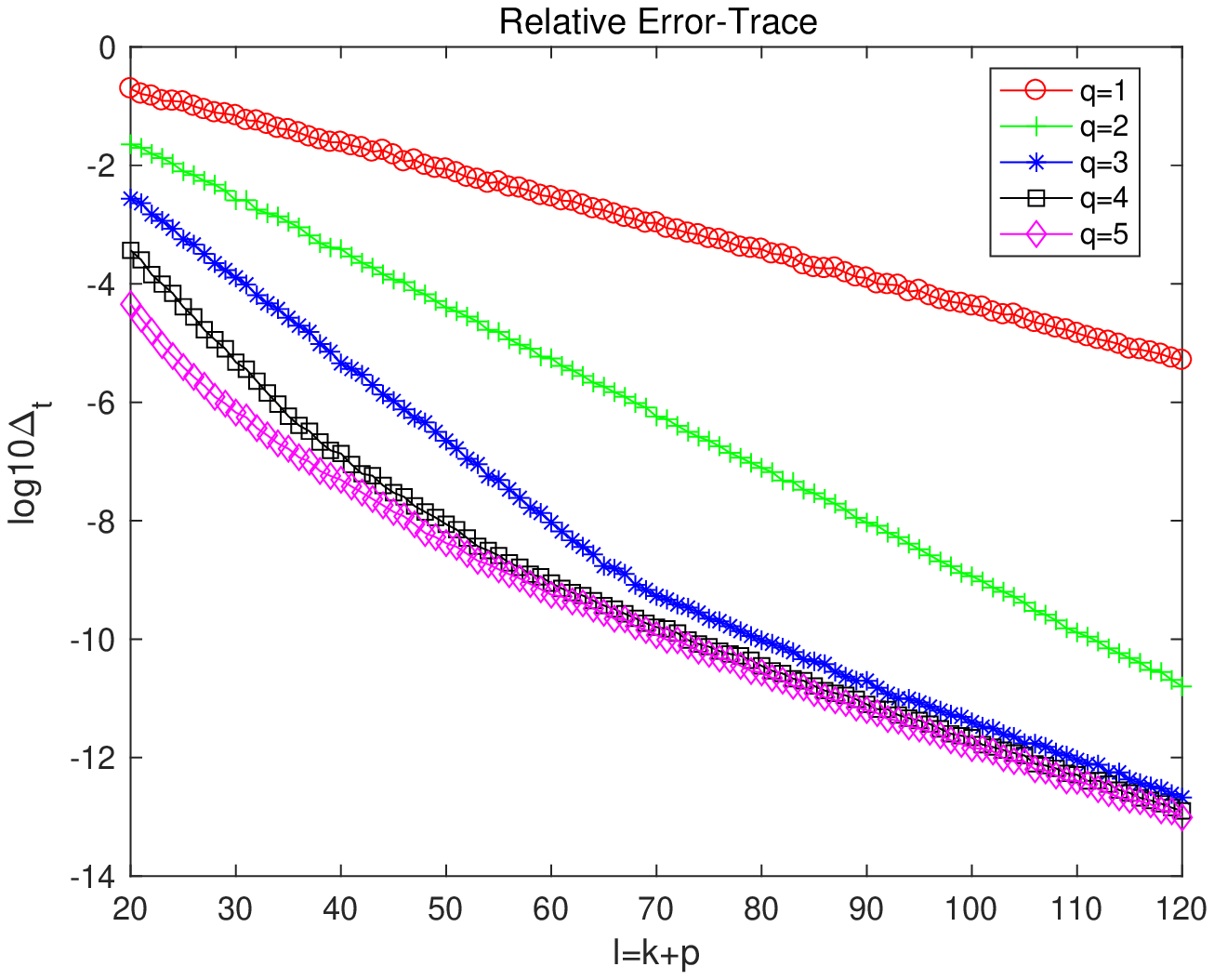}}
	\subfigure{%[SubfigureCaption]{
		%\label{Fig.sub.2}
		\includegraphics[height=48mm,width=55mm]{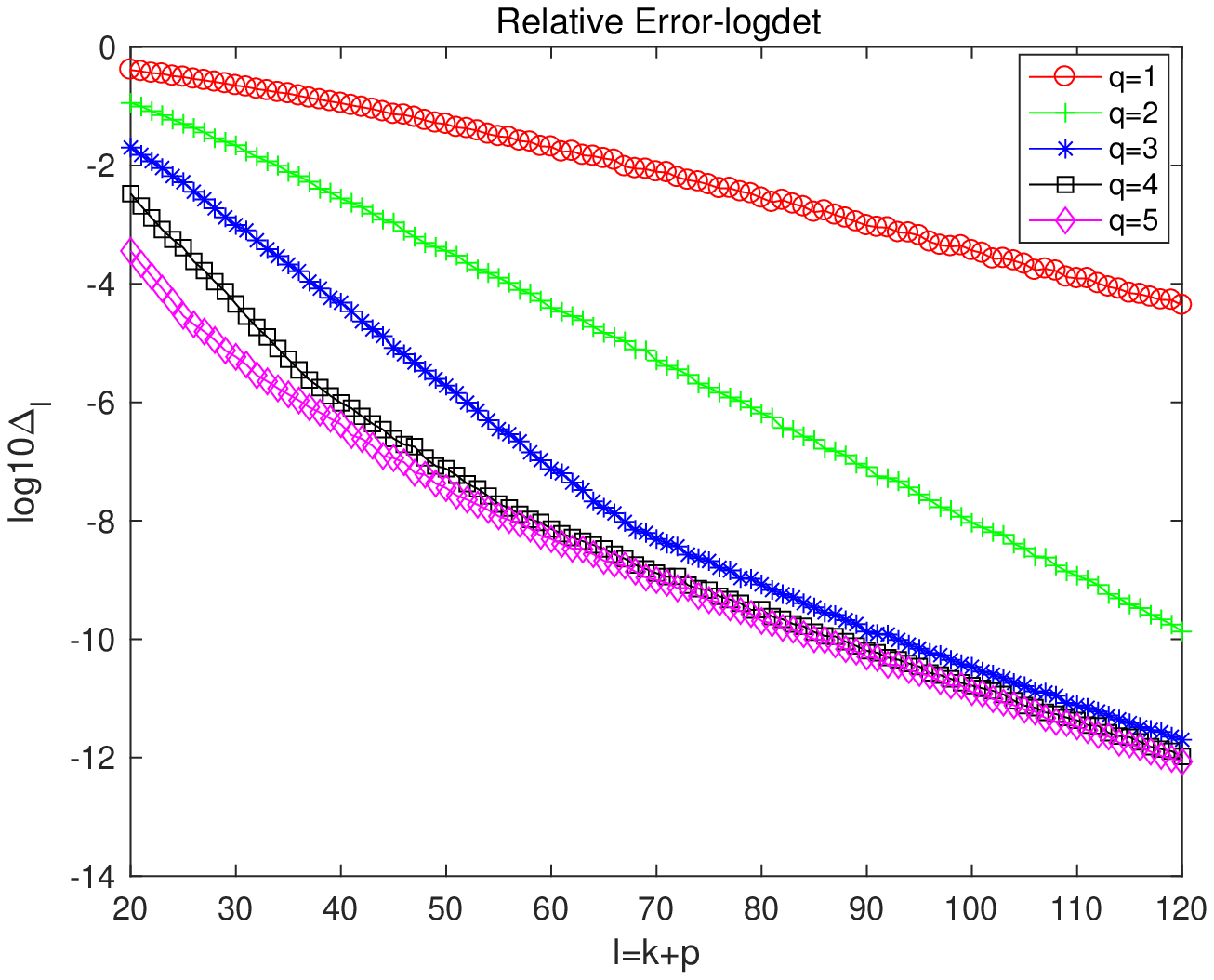}}
	\caption{%Accuracy of proposed estimators on a matrix with geometrically decaying eigenvalues.
		Accuracy of (left) trace and (right) log-determinant estimators produced by Algorithm \ref{Al2} for small matrix with $q$ varying from $1$ to $5$. The relative error is plotted against the sample size }
	\label{Fig2}
\end{figure}

The second experiment is used to test the effect of block Krylov space iteration parameter $q$ on the algorithms. Numerical results are displayed in Fig. \ref{Fig2}, which shows that the accuracy of both the trace and log-determinant estimators increases as the parameter $q$ increases for a fixed target rank $k$.  However, the growth is slowing as $q$ is increasing. As pointed out in \cite{saibaba}, this is because the overall error is dominated by ${\rm Tr}(\Lambda_2)$ and $\log\det(I_{n-k}+\Lambda_2)$.

\begin{figure}[H]
	\centering
	\subfigure{%[SubfigureCaption]{
		%\label{Fig.sub.1}
		\includegraphics[height=48mm,width=55mm]{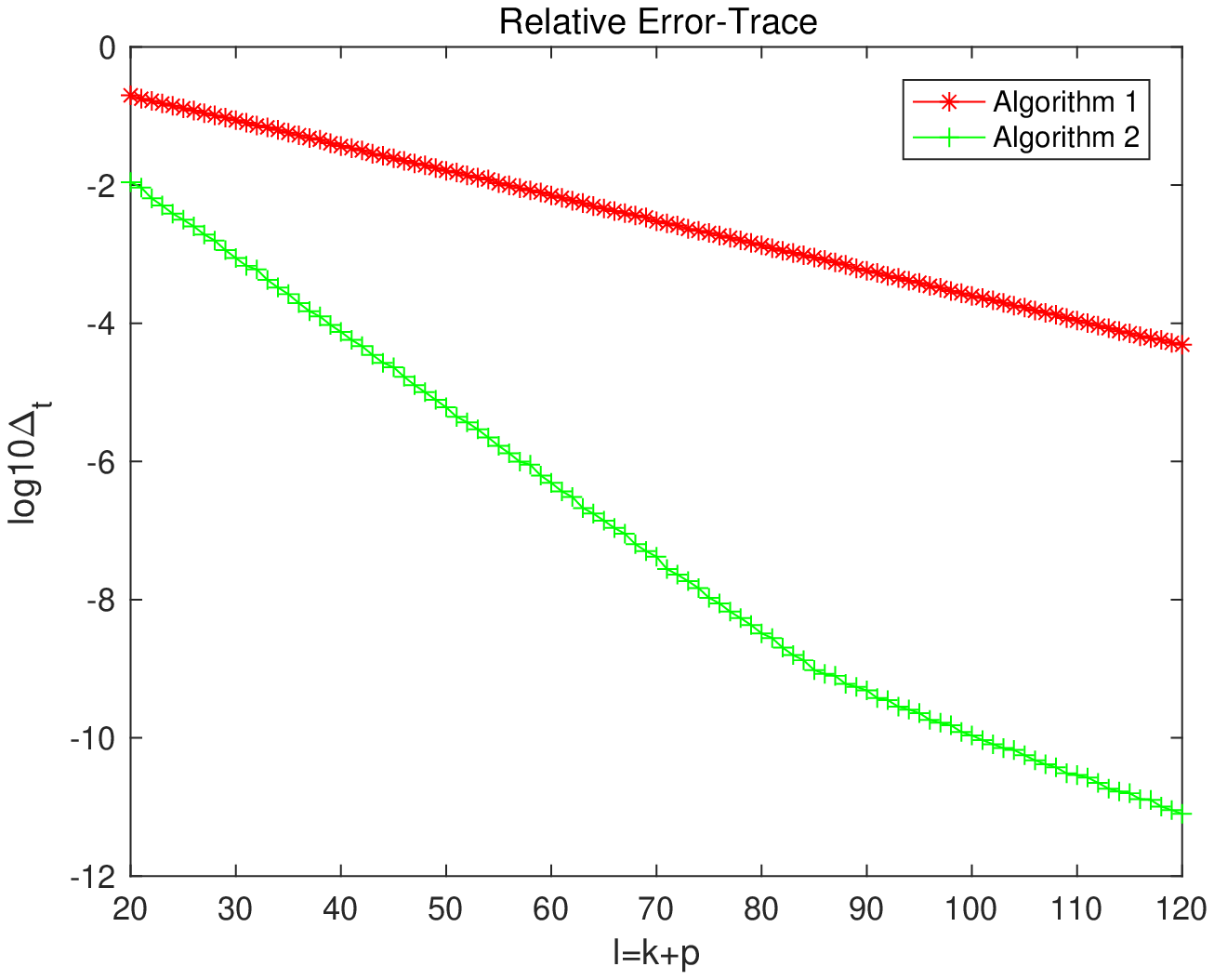}}
	\subfigure{%[SubfigureCaption]{
		%\label{Fig.sub.2}
		\includegraphics[height=48mm,width=55mm]{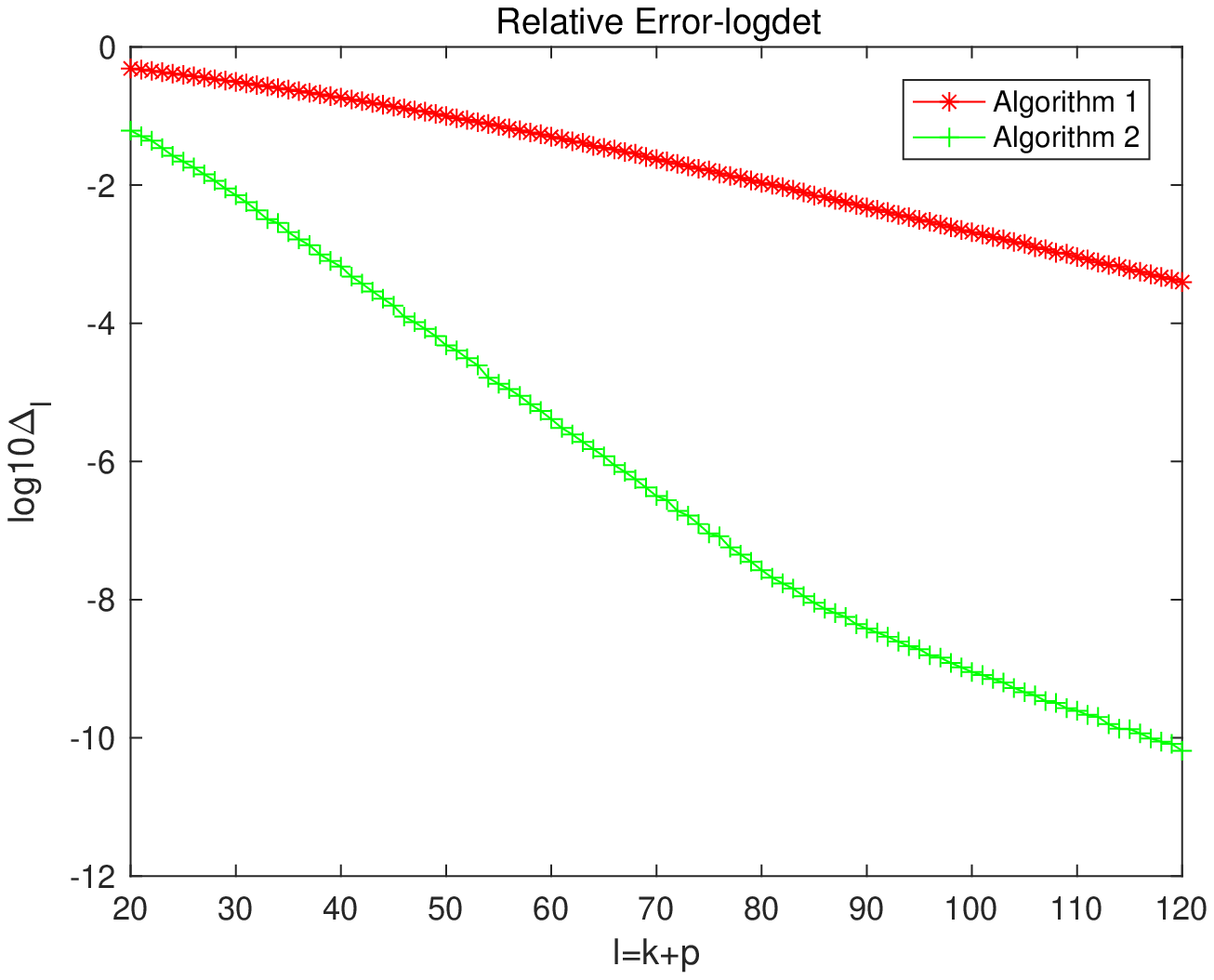}}
	\caption{%Accuracy of proposed estimators on a matrix with geometrically decaying eigenvalues.
		Comparisons of  Algorithms \ref{Al1} and \ref{Al2} for (left) trace and (right) log-determinant estimators for small matrix. The relative error is plotted against the sample size}
	\label{Fig3}
\end{figure}

In the third experiment, we compare Algorithms \ref{Al1} and \ref{Al2} for a special setting on sampling parameter $p$, block Krylov space iteration parameter $q$, and eigenvalue gap $\tau$. Numerical results are displayed in Fig. \ref{Fig3}, from which we can find that the estimators produced by Algorithm \ref{Al2} is always more accurate than the corresponding ones produced by Algorithm \ref{Al1}, and the differences increase as the sample size increases. In this experiment, we set $q=3$. We also do the experiments for $q\geqslant 2$. The results are similar and the differences for fixed  target rank $k$ will increase when $q$ increases. Of course, the two algorithms behave the same when $q=1$.

\begin{figure}[H]
	\centering
	\subfigure{%[SubfigureCaption]{
		%\label{Fig.sub.1}
		\includegraphics[height=48mm,width=55mm]{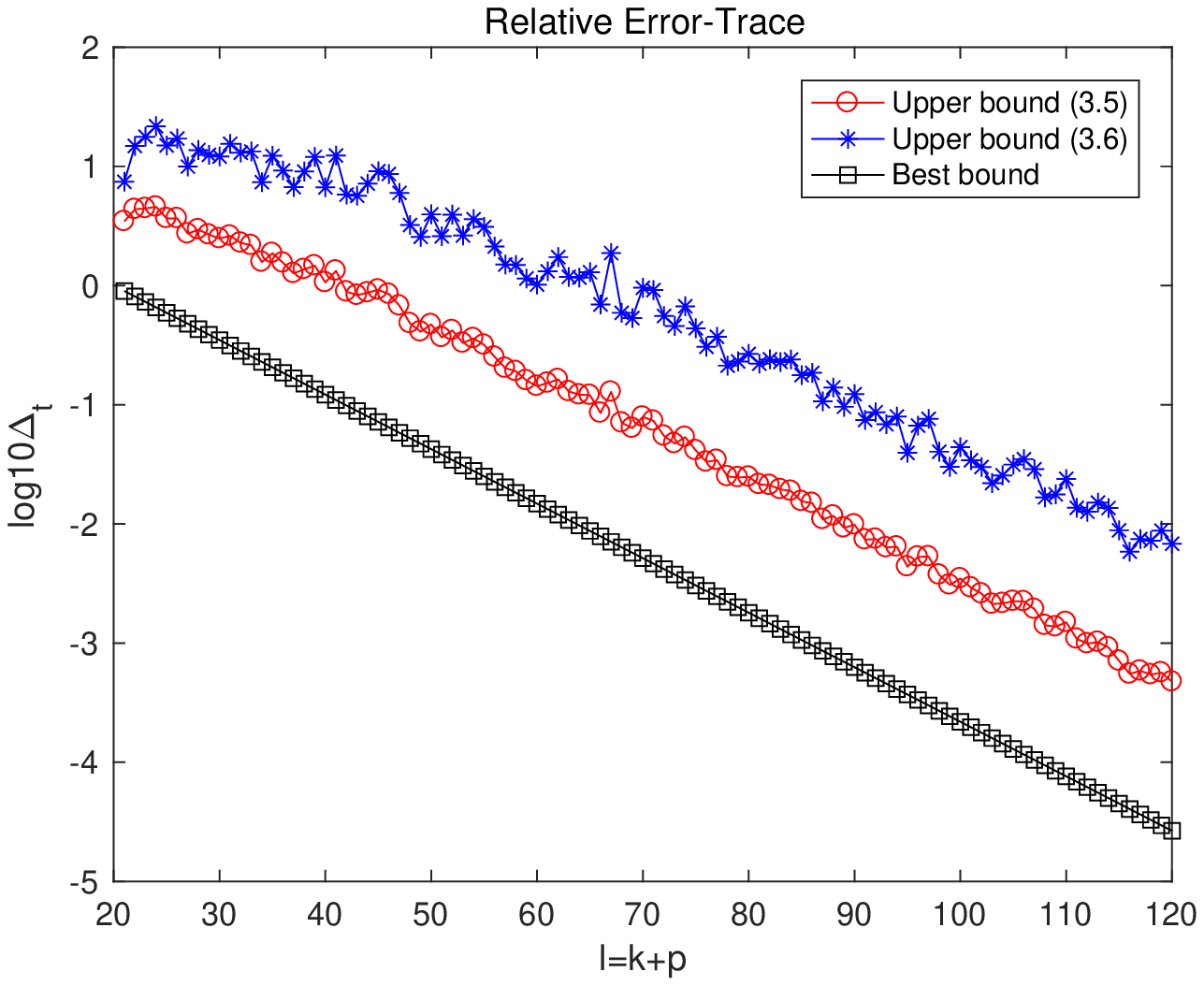}}
	\subfigure{%[SubfigureCaption]{
		%\label{Fig.sub.2}
		\includegraphics[height=48mm,width=55mm]{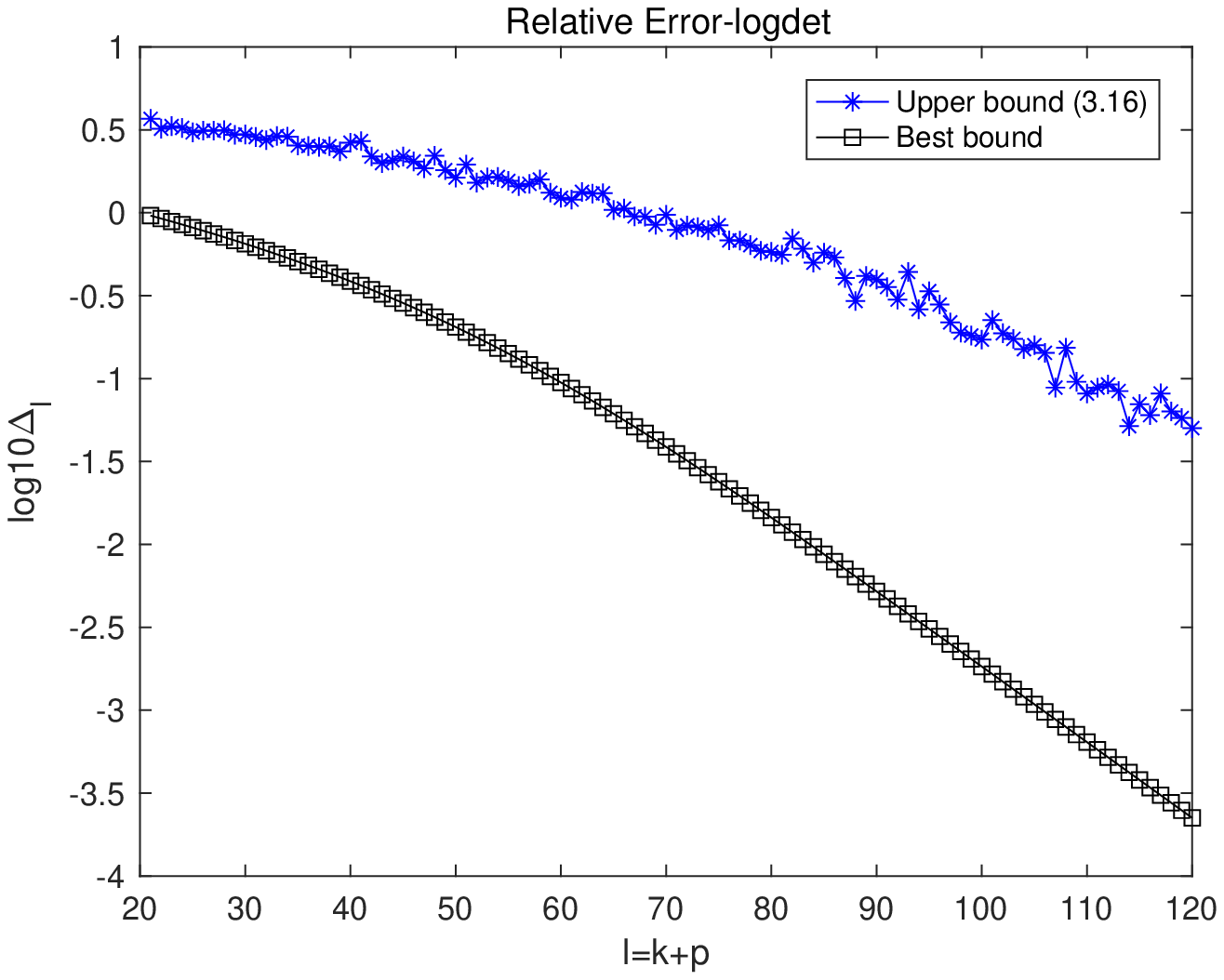}}
	\caption{%Accuracy of proposed estimators on a matrix with geometrically decaying eigenvalues.
		Comparisons of new error bounds and the best error bounds for  (left) trace and (right) log-determinant estimators for small matrix. The relative error is plotted against the sample size }
	\label{Fig4}
\end{figure}

In the fourth experiment, we test the accuracy of our structural error bounds. Specifically, we compare the error bounds \eqref{10000}, \eqref{10001}, and \eqref{40000} with the corresponding best error bounds ${\rm Tr}(\Lambda_2)$ and $\log\det(I_{n-k}+\Lambda_2)$. It should be clarified that the results displayed in Fig. \ref{Fig4} are all relative error bounds. That is, they are divided by ${\rm Tr}(A)$ and $\log\det(I_{n}+A)$, respectively. These results suggest that our bounds are effective. Especially, the bounds for trace estimator are qualitatively similar to the best error bound and the bound \eqref{10000} is also quantitatively within a factor of 10 of the best error bound.  %In addition, because, for the setting of this experiment, $\tau$ is close to 1, and
%$$\left \| \widehat{\Omega}_{2}\widehat{\Omega}_{1}^{\dagger}\right\|_{2}\geqslant\frac{\lambda _{k}}{\lambda _{k+1}}T_{q-1}(\frac{2\lambda_{k}-\lambda_{k+1}}{\lambda_{k+1}}),$$
%the bound \eqref{10000} is tighter than \eqref{10001}, and more accurate.
For log-determinant estimator, the differences between the error bound  \eqref{40000} and the best bound become a little larger when sample size increases.

\subsection{Medium sized matrices}\label{Sec5.2}
The test matrix $A$ is defined as
\begin{align}\label{503}
A\equiv \sum_{j=1}^{40}\frac{h}{j^{2}}x_{j}x_{j}^{T}+\sum_{j=41}^{300}\frac{l}{j^{2}}x_{j}x_{j}^{T}.
\end{align}
In contrast to \cite{saibaba}, we set the dimension of the sparse vectors $ x_{j}$ with random nonnegative entries to be $20000 $. So the matrix $A$ is of size $20000\times 20000$. Note that $  x_{j} $ are not orthonormal. They are generated by the Matlab command\\ \textsf{$x_{j}=sprand(20000,1,0.025)$}. By setting suitable values of $h$ and $l$, we do the following four specific numerical experiments.
%Furthermore, the eigenvalues decay like $ \frac{1}{j^{2}} $ with a gap at index 40, and its magnitude depends on the ratio $ \frac{h}{l} $. The
%exact rank of this matrix is 300. 
\begin{enumerate}
	\item Test the performance of Algorithm \ref{Al2} when $p = 20$, $h = 10$, $l= 1$, and $q$ varies from 1 to 5.
	\item Test the performance of Algorithm \ref{Al2} when $p = 20$, $h = 1000$, $l= 1$, and $q$ varies from 1 to 5.
	\item Compare Algorithms \ref{Al1} and \ref{Al2}, when $p = 20$, $q = 3$, $h = 10$, and $l = 1$.
	\item Compare Algorithms \ref{Al1} and \ref{Al2}, when $p = 20$, $q = 3$, $h = 1000$, and $l = 1$. %, and compare the
\end{enumerate}

\begin{figure}[H]
	\centering
	\subfigure{%[SubfigureCaption]{
		%\label{Fig.sub.1}
		\includegraphics[height=48mm,width=55mm]{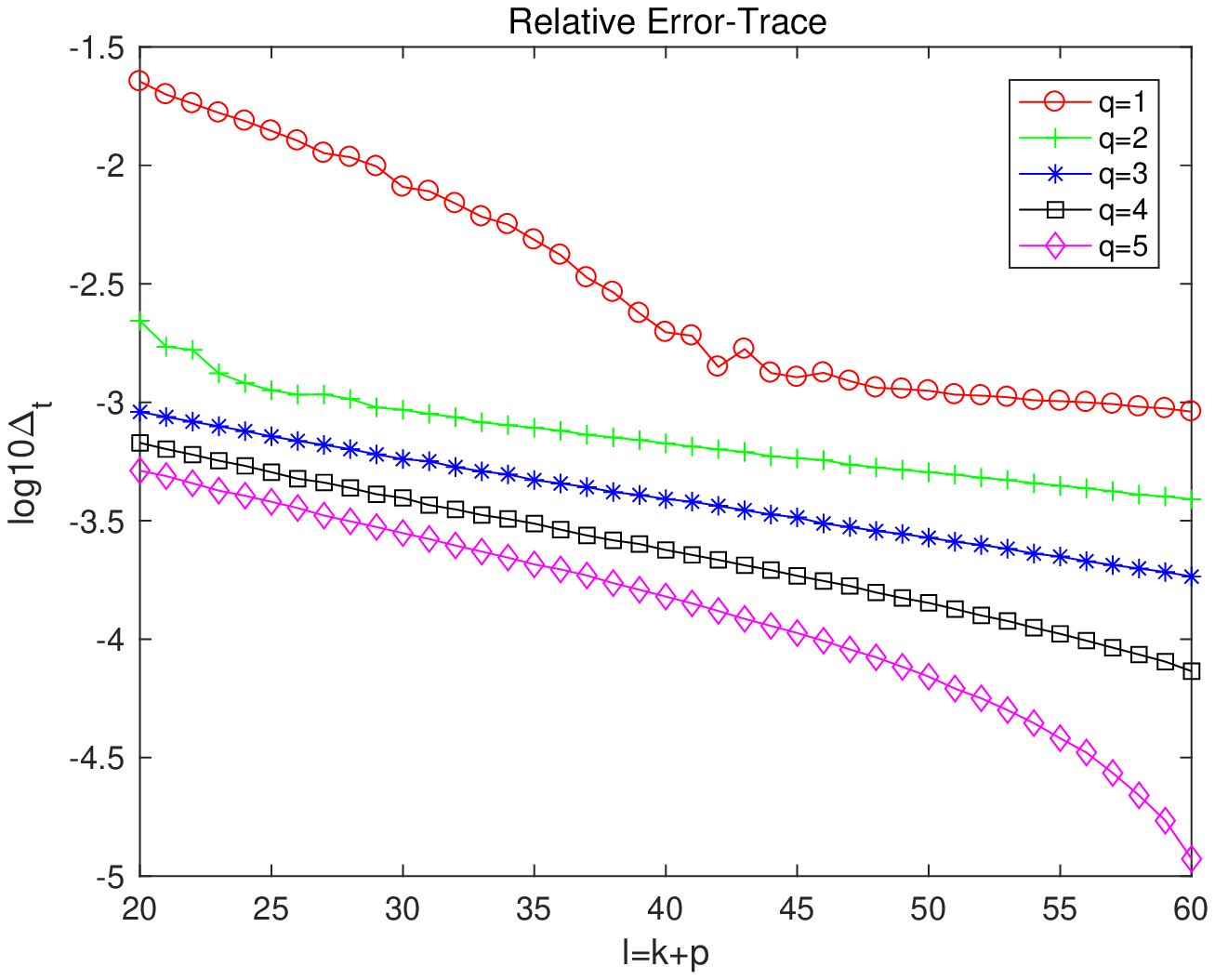}}
	\subfigure{%[SubfigureCaption]{
		%\label{Fig.sub.2}
		\includegraphics[height=48mm,width=55mm]{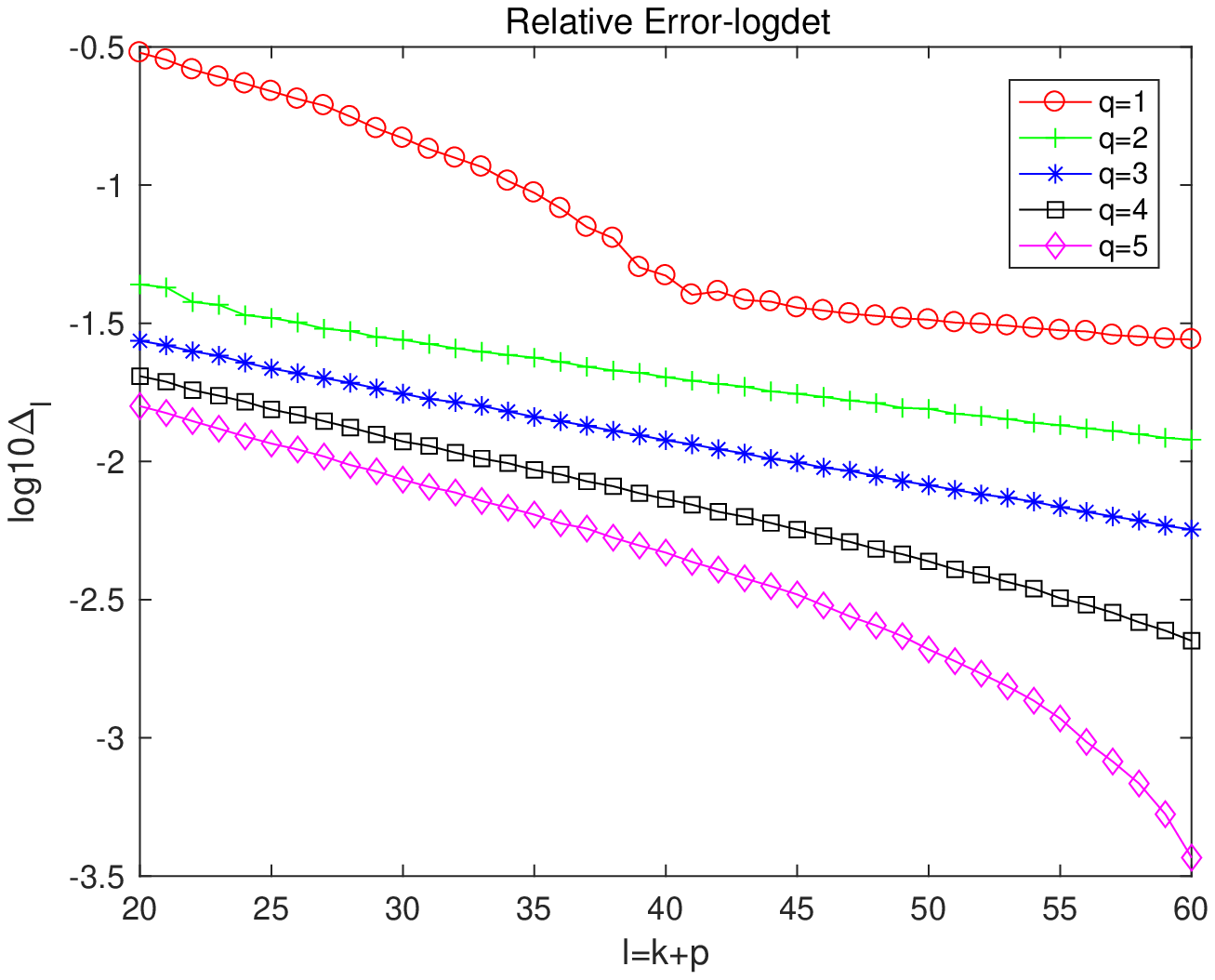}}
	\caption{%Accuracy of proposed estimators on a matrix with geometrically decaying eigenvalues.
		Accuracy of  (left) trace and (right) log-determinant estimators for medium matrix with $h = 10$ and $q$ varying  from $1$ to $5$. The relative error is plotted against the sample size}
	\label{Fig5}
\end{figure}

\begin{figure}[H]
	\centering
	\subfigure{%[SubfigureCaption]{
		%\label{Fig.sub.1}
		\includegraphics[height=48mm,width=55mm]{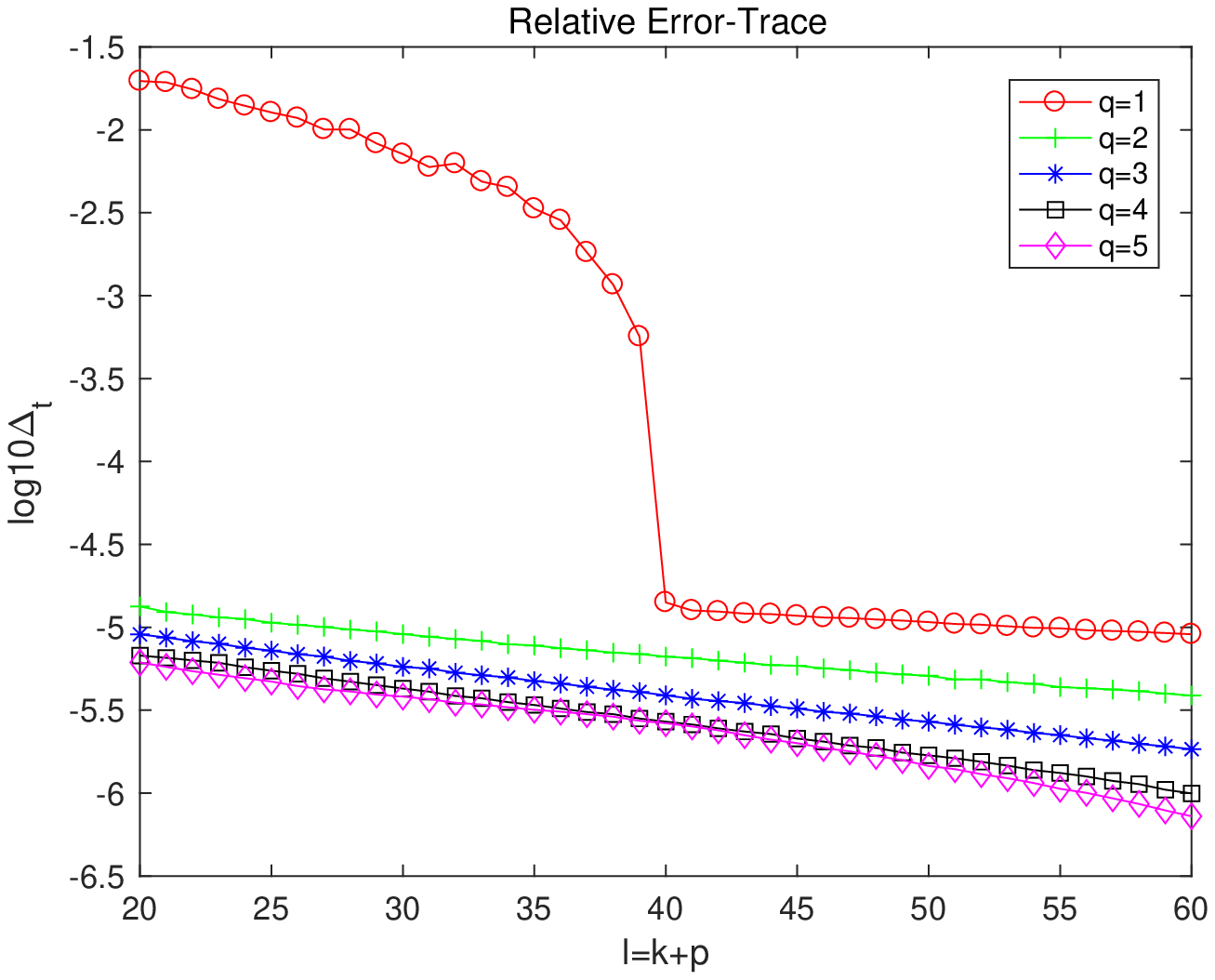}}
	\subfigure{%[SubfigureCaption]{
		%\label{Fig.sub.2}
		\includegraphics[height=48mm,width=55mm]{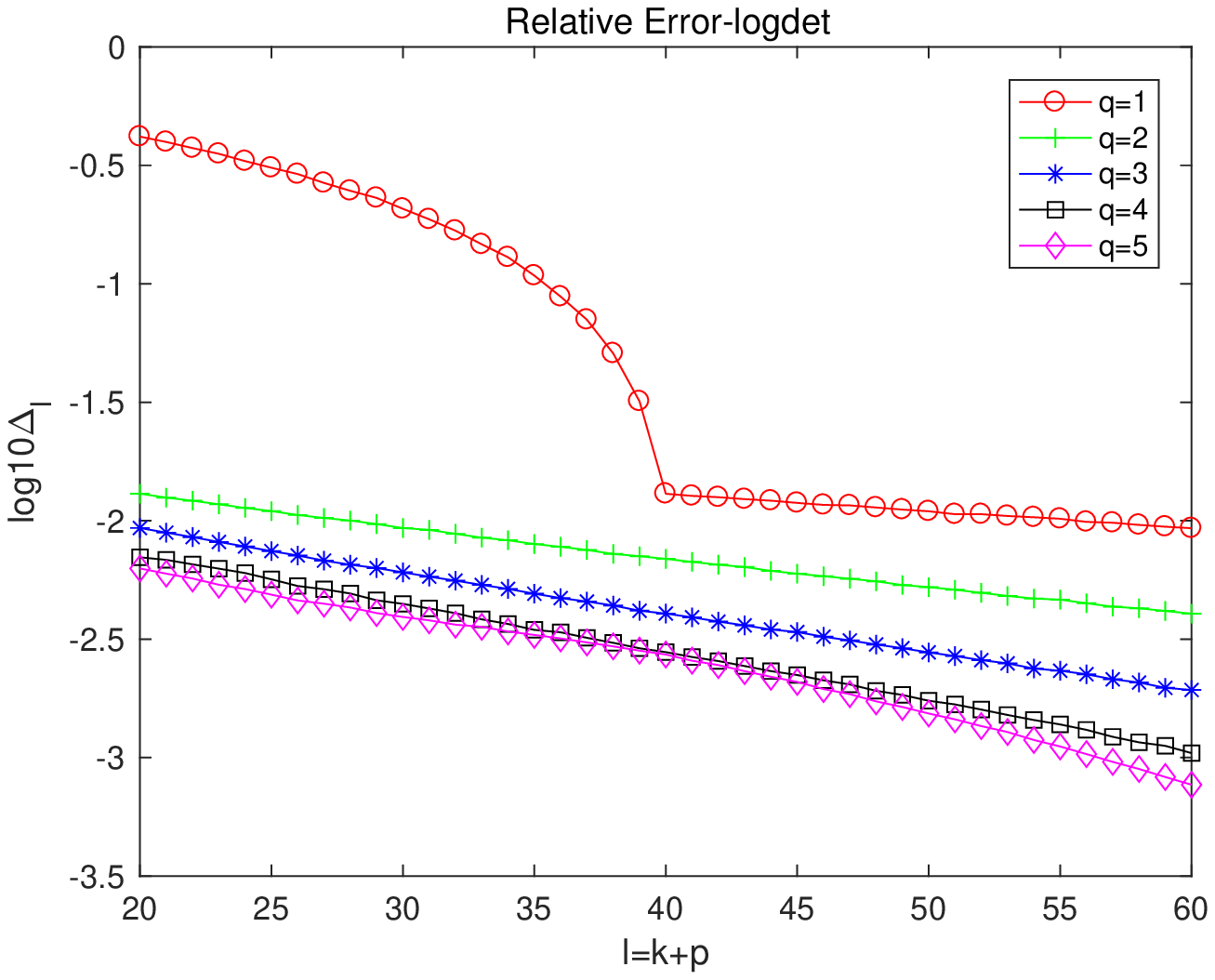}}
	\caption{%Accuracy of proposed estimators on a matrix with geometrically decaying eigenvalues.
		Accuracy of  (left) trace and (right) log-determinant estimators for medium matrix with $h = 1000$ and $q$ varying  from $1$ to $5$.  The relative error is plotted against the sample size}
	\label{Fig6}
\end{figure}

\begin{figure}[H]
	\centering
	\subfigure{%[SubfigureCaption]{
		%\label{Fig.sub.1}
		\includegraphics[height=48mm,width=55mm]{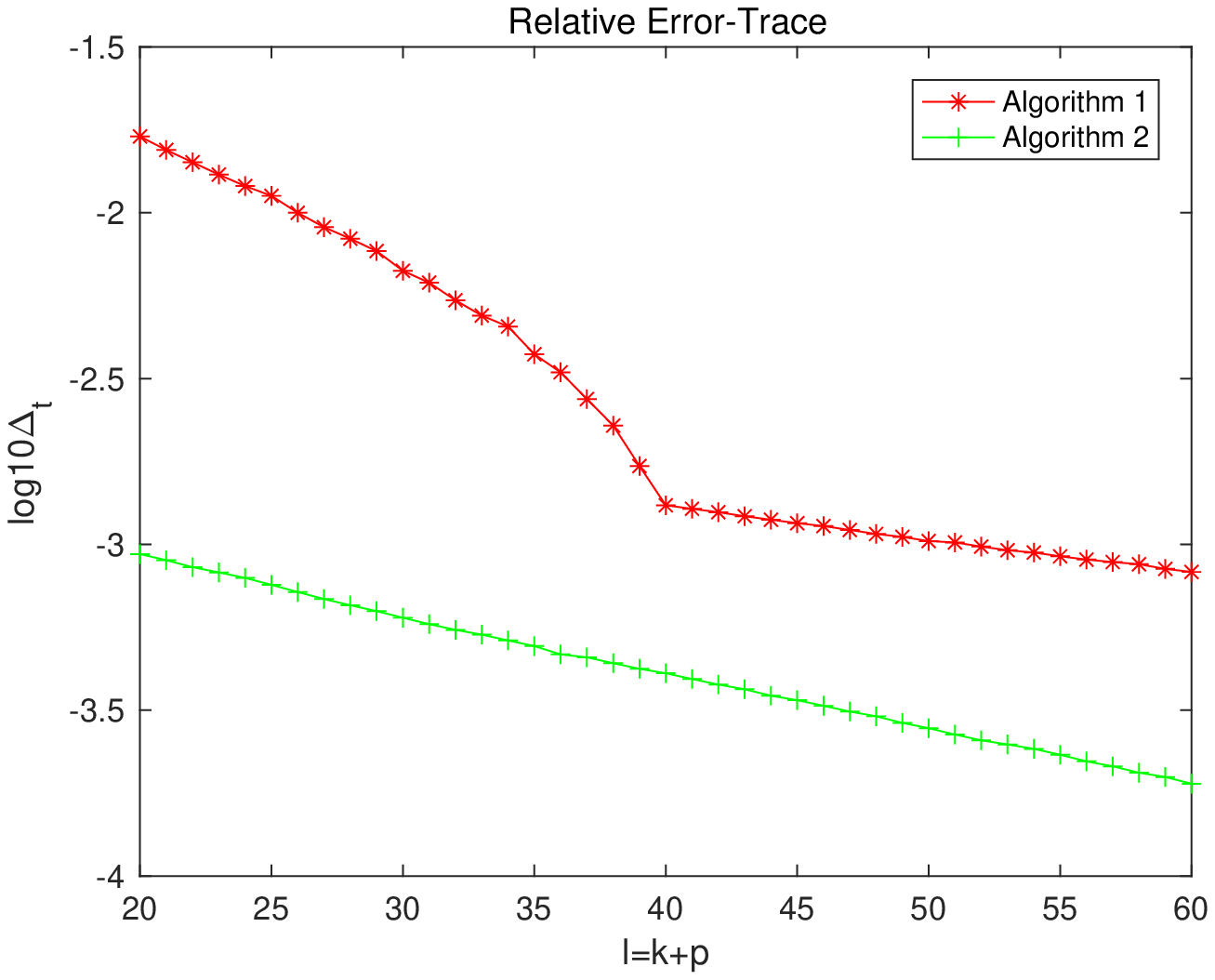}}
	\subfigure{%[SubfigureCaption]{
		%\label{Fig.sub.2}
		\includegraphics[height=48mm,width=55mm]{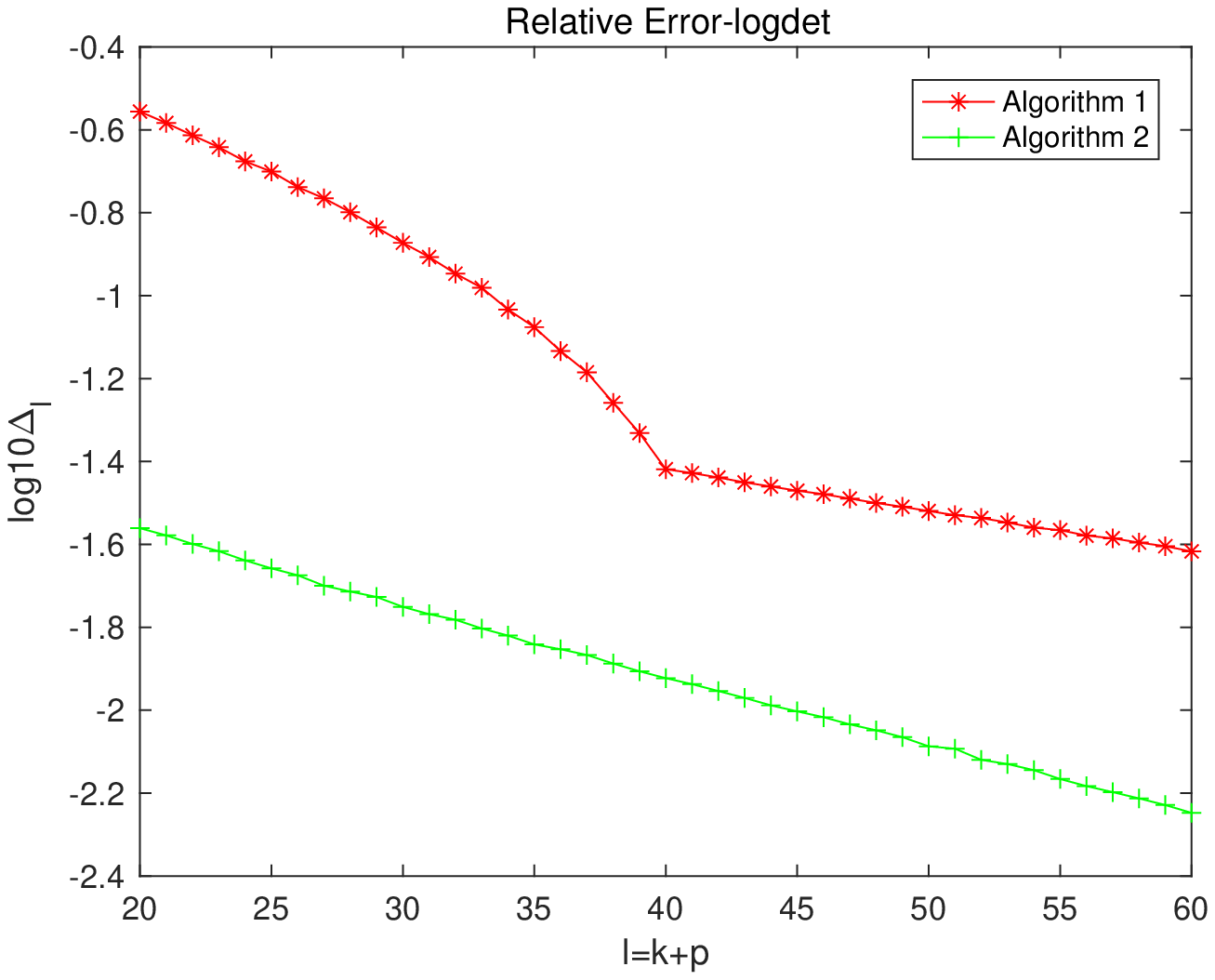}}
	\caption{%Accuracy of proposed estimators on a matrix with geometrically decaying eigenvalues.
		Comparisons of  Algorithms \ref{Al1} and \ref{Al2} for  (left) trace and (right) log-determinant estimators for medium matrix with $h=10$. The relative error is plotted against the sample size  }
	\label{Fig7}
\end{figure}

\begin{figure}[H]
	\centering
	\subfigure{%[SubfigureCaption]{
		%\label{Fig.sub.1}
		\includegraphics[height=48mm,width=55mm]{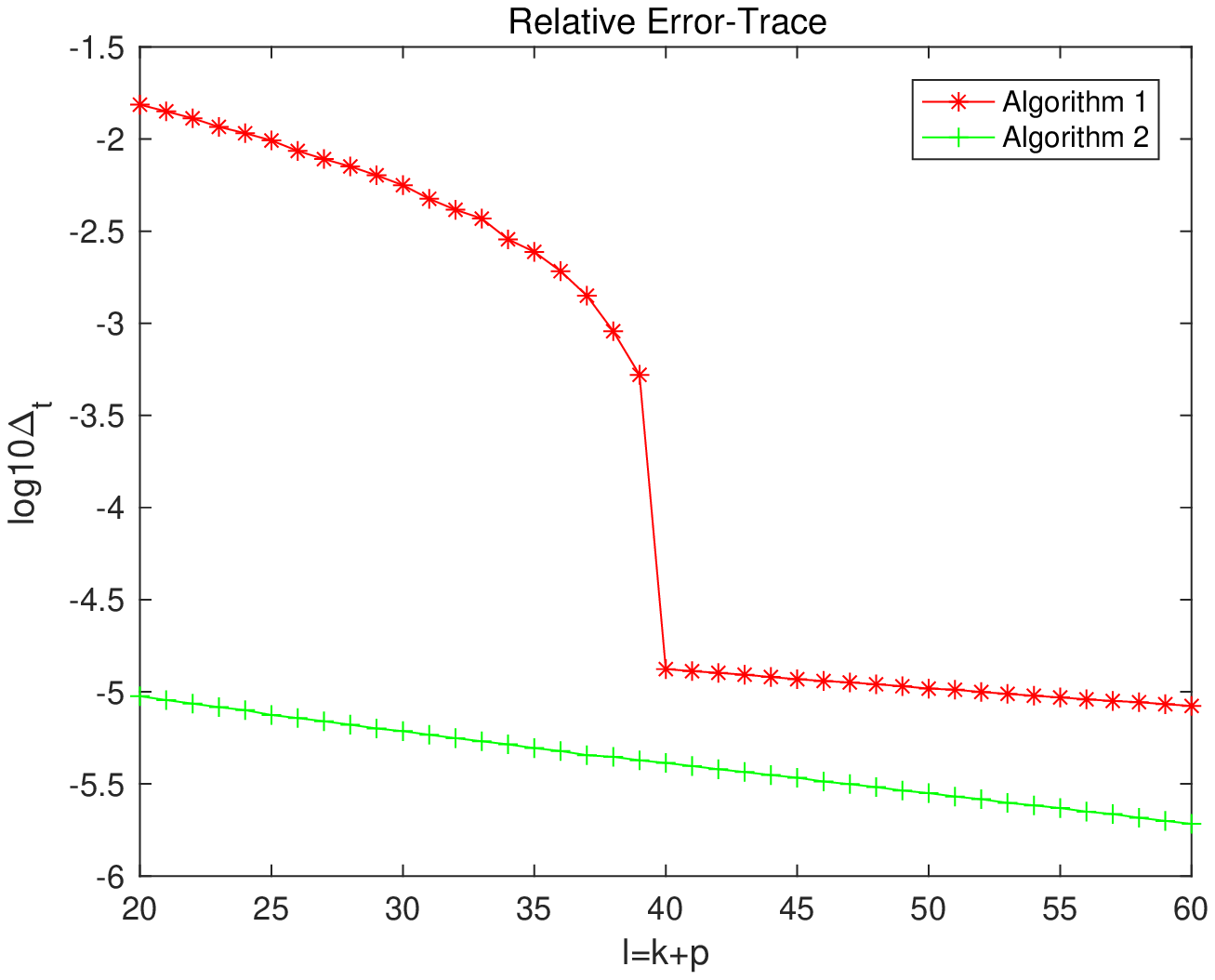}}
	\subfigure{%[SubfigureCaption]{
		%\label{Fig.sub.2}
		\includegraphics[height=48mm,width=55mm]{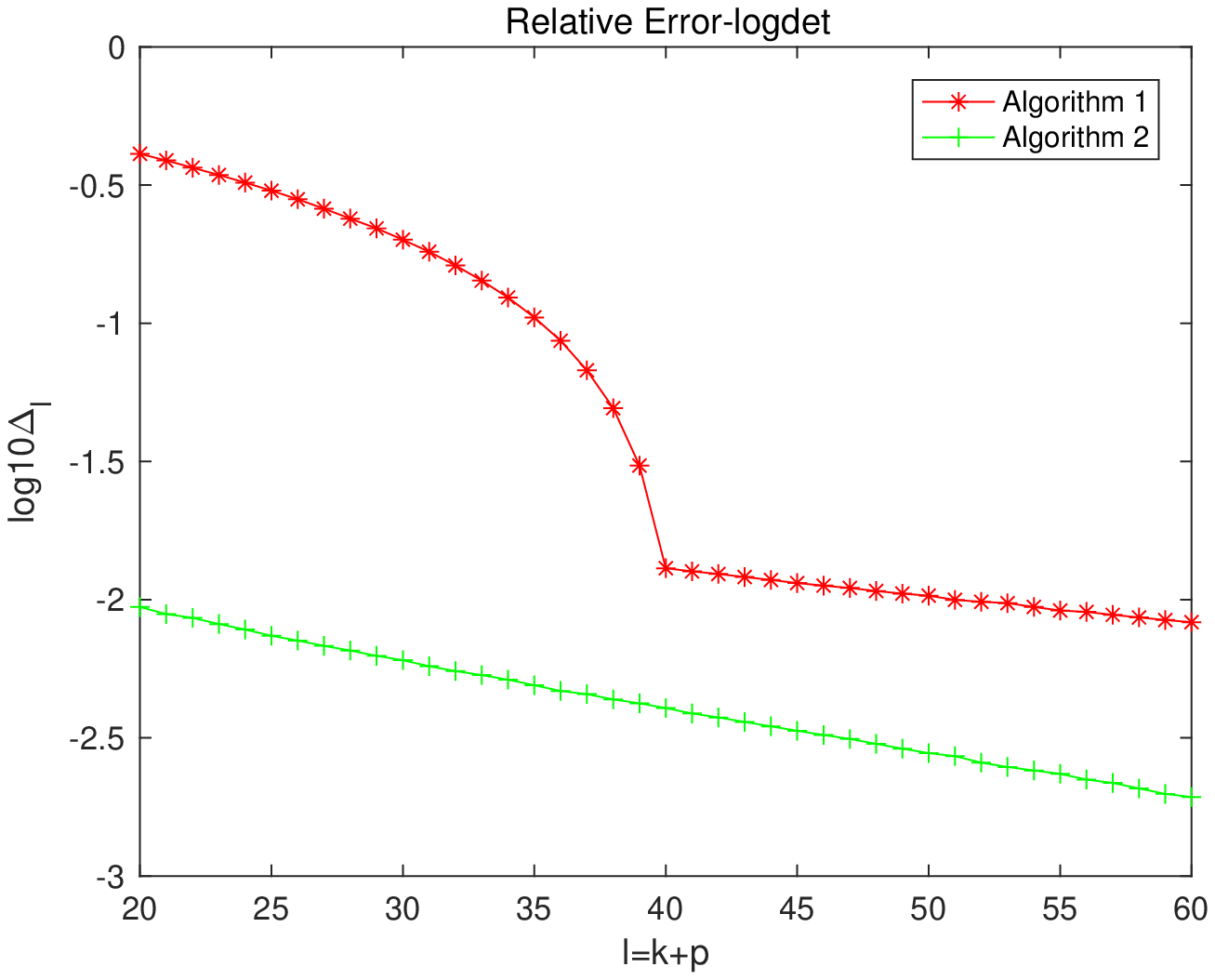}}
	\caption{%Accuracy of proposed estimators on a matrix with geometrically decaying eigenvalues.
		Comparisons of  Algorithms \ref{Al1} and \ref{Al2} for  (left) trace and (right) log-determinant estimators for medium matrix with $h=1000$. The relative error is plotted against the sample size }
	\label{Fig8}
\end{figure}

Numerical results of these experiments are displayed in Figs. \ref{Fig5}--\ref{Fig8}, respectively, which show the similar results found in the experiments in Section \ref{Sec5.1}. That is, the accuracy of both estimators increases as the parameter $q$ increases for a fixed target rank $k$,  the growth is slowing as $q$ is increasing, and the estimators produced by Algorithm \ref{Al2} is always more accurate than the corresponding ones produced by Algorithm \ref{Al1}. As pointed out in \cite{saibaba}, the accuracy of Algorithm \ref{Al1} improves considerably around the location of eigenvalue jump, and the larger the jump, the greater the improvement. In contrast, the accuracy of our algorithm is quite high in all the locations. This is mainly because we use the information discarded by Algorithm \ref{Al1}. These information improves the accuracy of estimators greatly.

\section{Concluding remarks}
In this paper, we present new randomized algorithms for estimating trace and log-determinant of Hermitian positive semi-definite matrices defined implicitly. Numerical experiments show that the performance of the new algorithms is better than that for algorithms given in \cite{saibaba}. We also provide rigorous error bounds for our trace and log-determinant estimators. They are tighter than the corresponding ones from \cite{saibaba} in most of cases. %However, unlike \cite{saibaba}, we don't consider the Rademacher random matrices as the starting guess of our algorithms. 

To achieve the information from $\mathcal{K}_q$, we adopt a popular method \cite{drineas}, that is, we  consider an element of $\mathcal{K}_q$: $\phi\left(A\right)\Omega$. The method has a drawback that it requires $l\geqslant k$. The requirement can be relaxed by considering a method from \cite{yuan}. However, it will be difficult to investigate the expectation and concentration error bounds of estimators for this method. In addition, for concise and comparing with the error bounds in \cite{saibaba} directly, we partition $\Lambda$ in \eqref{201} into two blocks. As done in \cite{gu2015,yuan}, we can partition $\Lambda$ into three blocks or four blocks. That is, we make an artificial gap and consider a cluster of some eigenvalues. The error bounds for this setting will be more flexible compared with the results obtained in this paper. We will consider these problems in the future work.

%\begin{acknowledgements}
%\end{acknowledgements}




\end{document}